\documentclass[11pt,amstex]{article}
\usepackage{amsmath}
\usepackage{amsfonts}
\usepackage{amscd}
\usepackage{amsthm}
\usepackage[dvips]{color}

\input xypic
\xyoption{all}

\def\a{\alpha}
\def\b{\beta}

\def\lam{\lambda}
\def\si{\sigma}

\def\be{\begin{equation}}
\def\ee{\end{equation}}
\def\bear{\begin{eqnarray}}
\def\eear{\end{eqnarray}}
\def\best{\begin{eqnarray*}}
\def\eest{\end{eqnarray*}}
\def\pf{{\bf Proof}: }

\newtheorem{theorem}{Theorem}[section]
\newtheorem{prop}[theorem]{Proposition}

\newtheorem{lemma}[theorem]{Lemma}

\newtheorem{defn}[theorem]{Definition}

\newtheorem*{thmA}{Theorem A}
\newtheorem*{thmB}{Theorem B}

\newtheorem{remark}[theorem]{Remark}
\newenvironment{rem}{\begin{remark}\rm}{\end{remark}}
\newtheorem{example}[theorem]{Example}

\def\non{\noindent}
\def\pf{\non {\bf Proof. }}
\def\qed{\nopagebreak \hskip .1in { $\Box$ }\penalty10000 %
\hskip\parfillskip \par  }
\def\ra{\rightarrow}

\def\del{\overline \partial}
\def\bd{\partial}

\def\cok{\mbox{coker\,}}

\def\P{{\mathbb P}}
\def\Q{{\mathbb Q}}
\def\cx{{\mathbb C}}

\def\im{\mbox{\rm Im}}

\def\Hom{{\rm Hom}}

\def\J{{\cal J}}

\def\O{{\cal O}}

\def\K{{\cal K}}

\def\M{{\cal M}}

\newcommand{\CM}{\overline{{\cal M}}}
\newcommand{\CU}{\overline{{\cal U}}}

\def\Map{{\cal{M}}{\it{ap}}}

\def\F1{{\mathbb F}_1}
\def\la{\langle}
\def\ra{\rangle}

\def\F{\mathbb F}

\def\vir{\scriptstyle vir}

\def\sv{\scriptscriptstyle V}

\def\si{\scriptscriptstyle I}

\def\s0{{\scriptscriptstyle 0}}

\def\TM{\widetilde{\M}}

\def\ve{\varepsilon}

\begin{document}

\title{\bf Sum Formulas for Local Gromov-Witten Invariants of Spin Curves}
\vskip.2in

\author{ Junho Lee \\ University of Central Florida \\ Orlando, FL 32816
}

\date{\empty}

\maketitle

\begin{abstract}
\medskip

Holomorphic 2-forms on K\"{a}hler surfaces  lead to ``Local Gromov-Witten invariants'' of spin curves.
This paper shows how to derive sum formulas for such local GW invariants from
the sum formula for GW invariants of certain ruled surfaces. These sum formulas
also verify the Maulik-Pandharipande formulas that were recently proved by Kiem and Li.

\end{abstract}

\vspace{1cm}

\setcounter{section}{0}

Let $X$ be a K\"{a}hler surface with a holomorphic 2-form $\a$. The real part of  $\a$, also denoted by $\a$,
 then induces an almost complex structure on $X$\,:
\begin{equation}\label{Ja}
J_\a\,=\,(Id+JK_\a)^{-1}J(Id+JK_\a).
\end{equation}
Here $J$ is the K\"{a}hler structure on $X$ and $K_\a$ is the endomorphism of $TX$ defined by
the formula $\la u,K_\a v\ra =\a(u,v)$ where $\la\,,\,\ra$ is the K\"{a}hler metric on $X$.
The almost complex structure $J_\a$ satisfies a remarkable {\em Image Localization Property\,:}
\begin{itemize}
\item {\em
if $f$ is a $J_\a$-holomorphic map into $X$  that represents a non-zero (1,1) class  then
the image of $f$ lies in the zero set $D$ of the holomorphic 2-form $\a$}
\end{itemize}
For simplicity, assume $X$ is a (minimal) surface of general type and $D$ is smooth.
The normal bundle $N$ to $D$ is then a theta characteristic on $D$, i.e., $N$ is a square root of the
canonical bundle of $D$.
The pair $(D,N)$ is called a  {\em spin curve} of genus $h$ where $h$ is the genus of $D$.
The total space $N_D$ of $N$ has a tautological holomorphic 2-form $\a$ that induces,
by the same manner as in (\ref{Ja}), an almost complex structure $J_\a$ on
$N_D$ satisfying the image localization property, namely
\begin{equation*}\label{ILL-N}
\CM_{g,n}(N_D,d[D],J_\a)\,=\,\CM_{g,n}(D,d).
\end{equation*}
Consequently, the moduli space of $J_\a$-holomorphic maps is compact, so it represents a (virtual) fundamental class
that defines local GW invariants  of the spin curve $(D,N)$.
These local GW invariants depend only on the genus $h$ and the parity $p\equiv h^0(N)$ (mod 2)
and GW invariants of K\"{a}hler surfaces
with $p_g>0$ are the sum of  local GW invariants  associated to spin curves \cite{LP1}.

GW invariants count maps from connected domains, while Gromov-Taubes invariants
count maps from not necessarily connected domains.
These GT invariants can be obtained from GW invariants.
Maulik and Pandharipande \cite{MP} gave fascinating conjectural formulas for (descendent) local GT invariants
of spin curve $(D,N)$ for low degrees ($d=1,2$)\,:
\begin{equation}\label{MPf}
\begin{array}{l}
{\displaystyle
GT_{1}^{loc,h,p}\big(\,\prod_{i=1}^n \tau_{k_i}(F^*)\,\big)\,=\,
(-1)^p\,\prod_{i=1}^n
\,\frac{k_i!}{(2k_i+1)!}(-2)^{-k_i}
\phantom{a_{\displaystyle \int}} } \\
{\displaystyle
GT_{2}^{loc,h,p}\big(\,\prod_{i=1}^n \tau_{k_i}(F^*)\,\big)\,=\,
(-1)^p\, 2^{h+n-1}\, \prod_{i=1}^n\,
\frac{k_i!}{(2k_i+1)!}(-2)^{k_i}
}
\end{array}
\end{equation}
(see Section~\ref{S:LGW} for definition of  descendant local GT invariants).
Kiem and Li \cite{KL1,KL2,KL3} have since proved these formulas using their algebro-geometric version of local invariants.
Observing the formulas (\ref{MPf}) for genus zero spin curve directly follow from Proposition 2 of \cite{FP},
they showed the following reduction theorem\,:
\begin{itemize}
\item
{\em
for low degrees ($d=1,2$)
local GT invariants of higher genus spin curves can be reduced to local GT invariants of genus zero spin curve.
}
\end{itemize}
Their proof uses a sum formula (Theorem 4.2 of \cite{KL1}) for degeneration obtained by   certain blow-up
plus explicit calculation of the invariant  $GT_{2}^{loc,h,p}(\tau(F^*))$.

\medskip
\bigskip
The aim of this paper is twofold.
First, we give a  proof of the sum formula  for
degeneration  by blow-up in the context of symplectic geometry --- this sum formula is the same as Kiem and Li's sum formula except
for the constraints of relative invariants (see Theorem~A below).
Second,  we also prove new sum formulas  for
degeneration of  spin curves in the case of degree $d=2$ (see Theorem~B below).
Then, in the proof of the above reduction theorem (cf. Section 4 of \cite{KL1}),
we can replace  the calculation of invariant $GT_2^{loc,h,p}(\tau)$ by
Theorem~B (see Section~\ref{S:MP}).

The novelty of our approach is to use GW invariants of ruled surfaces.
Unlike the algebro-geometric approach, our local GW
invariants of $(D,N)$ are, in fact, {\em local contributions} to
GW invariants of the ruled surface  $\P_h=\P(N\oplus \O_D)$ that count
maps whose images are close to the zero section $D$ of
$\P_h$. A small neighborhood $U$ of   $D$  in
$\P_h$ is isomorphic to some neighborhood of the zero
section in the total space $N_D$ of $N$. Together with this isomorphism and some
bump function, the tautological holomorphic 2-form on $N_D$ induces  an almost complex structure $J_\a$ on $\P_h$
satisfying the image localization property, namely
\begin{equation*}\label{ILL}
\CM_{g,n}(U,dS,J_\a)\,=\,\CM_{g,n}(D,d)
\end{equation*}
where $S$ is the section class of $\P_h$, i.e. $S=[D]$.
The moduli space of $J_\a$-holomorphic maps into $U$ thus
represents a (virtual) fundamental class that gives
the local GW invariants of $(D,N)$.
This description of local GW invariants
is well suited to easily adapt the arguments  in \cite{IP1,IP2}
to a version of sum formulas for  local GW invariants.
The relative local GW invariants are simply the local contributions to
the relative GW invariants of $\P_h$ that count maps into $U$ relative a fixed fiber of $\P_h$.
In terms of those relative invariants,
sum formulas for local invariants  directly
follow from the main argument of \cite{IP2} for some cases.

Our relative invariants are, however,  not given by (virtual) fundamental class
of relative moduli space that is needed to define descendent classes.
To get around this issue, we relate  descendent invariants to relative invariants with
$\phi_i$ classes that are the first Chern classes of the relative cotangent bundles
over the space of stable curves (see Proposition~\ref{dec=rel}). Then, we use those relative invariants to
show the sum formula  for degeneration  by blow-up.

\smallskip

Let $\F_0=\P^1\times E$ be a ruled surface over $E=\P^1$. Then there is a unique section of $\F_0$
that passes through a given point. This simple observation enable us to apply
the main argument of \cite{IP2} for the symplectic fiber sum
$\P_h=\P_h\#_{V}\F_0$ to obtain a sum formula for degeneration by blow-up.
For any partition $m=(m_1,\cdots,m_\ell)$,  we set
$$
\ell(m)\,=\,\ell,\ \ \ \ \ \ \ \ \ \
|m|\,=\,\prod\,m_i,
\ \ \ \ \ \ \ \ \ \ m!\,=\,|{\rm Aut}(m)|
$$
where ${\rm Aut}(m)$ is the symmetric group permuting equal parts of $m$.  In Section~\ref{S:Blow-Up}, we show\,:

\begin{thmA}
\label{Intro-A}
Let $d\ne 0$ and $n_1+n_2=n$. Then
\begin{align}\label{sum-bu}
&GT_{d}^{loc,h,p}\big(\,\prod_{i=1}^n \tau_{k_i}(F^*)\,\big) \notag \\
= \
&\frac{1}{(d!)^2}\,\sum_{m} \,\frac{|m|}{m!}\,
GT_{(1^d),m}^{loc,h,p}\big(\,\prod_{i=1}^{n_1} \phi_i^{k_i}(F^*)\,\big)
\cdot
GT_{m,(1^d)}^{\F_0}\big(\,\prod_{i=1}^{n_2} \phi_i^{k_{n_1+i}}(F^*)\,\big)
\end{align}
where the sum is over all partitions $m$ of $d$
(see Section~\ref{S:RLGW} for definition and notation of relative invariants).
\end{thmA}

\medskip
Let $\F_1=\P(\O(1)\oplus \O_E)$ be a ruled surface over $E$.
Unlike the case of $\F_0$, the infinite section plus a fiber represents the section  class represented
by  the zero section of $\F_1$.
This causes the main difficulty to derive general sum formulas of local GW invariants for degeneration of spin curves from
the symplectic fiber sums
$$
\P_h\,=\,\P_{h_1}\#_{V_1}\#\F_1\#_{V_2}\P_{h_2}\ \ \ \ \ \ \
\mbox{and}\ \ \ \ \ \ \
\P_h\,=\,\P_{h-1}\,\underset{V_1\sqcup\, V_2}{\#}\,\F_1.
$$
However, when degree $d=2$,  simple limiting arguments (see Section~\ref{S:Ruled-II} and Lemma~\ref{key-scd}) allow us to apply
the same arguments as in the proof of Theorem~A. In Section~\ref{S:SC-Sum}
and \ref{S:SC-Sum-2}, we show\,:

\begin{thmB}
\label{Intro-B}
\begin{itemize}
\item[]
\item[(a)] If $h=h_1+h_2$ and $p\equiv p_1+ p_2$ (mod 2) then we have
\begin{equation}\label{sum-deg-1}
GT_{(2)}^{loc,h,p}\,=\,
(-1)^{p_1}2^{h_1}\,  GT_{(2)}^{loc,h_2,p_2} \, +\,
(-1)^{p_2}2^{h_2}\, GT_{(2)}^{loc,h_1,p_1}\, -\,
(-1)^p 2^{h}\, GT_{(2)}^{loc,0,+}.
\end{equation}
\item[(b)] If $h\geq 2$ or if $(h,p)=(1,+)$ then we have
\begin{equation}\label{sum-deg-2}
GT_{(2)}^{loc,h,p}\,=\,
4\,GT_{(2)}^{loc,h-1,p}\,-\,(-1)^p 2^h\,GT_{(2)}^{loc,0,+}.
\end{equation}
\end{itemize}
\end{thmB}

\bigskip
\medskip
\non
{\bf Acknowledgments.} I am very thankful to Thomas H. Parker for valuable discussions
and to referees for corrections and useful comments. I also
thank Bumsig Kim and Young-Hoon Kiem for discussions and Rahul Pandharipande, Davesh Maulik
and Aleksey Zinger for their interest in this work.

\vskip 1cm

\setcounter{equation}{0}
\section{Moduli Spaces}
\label{S:Moduli}
\bigskip

This section introduces  moduli spaces of curves and maps.
For $2g+n\geq 3$, let $\CU_{g,n}\to \CM_{g,n}$ be the universal curve over the Deligne-Mumford space.
Lifting to the moduli space of Prym curves (\cite{Lo}, \cite{ACV}), one may assume that
 $\CM_{g,n}$ is a manifold  and
every connected $n$ marked nodal curve $C$ of (arithmetic) genus $g$ has a stabilization $st(C)\in \CM_{g,n}$
that is isomorphic to a fiber of $\CU_{g,n}$. After fixing an embedding $\CU_{g,n}\hookrightarrow \P^N$, one can  obtain
a map
\begin{equation*}
\phi: C\to st(C)\to \CU_{g,n}\to \P^N.
\end{equation*}
Let $(X,\omega)$ be a compact symplectic  manifold with
an $\omega$-tamed almost complex structure $J$.
A $C^1$-map $f:C\to X$ is stable if the energy
\begin{equation}\label{energy}
E(f,\phi)=\frac12\int |df|^2 + |d\phi|^2
\end{equation}
is positive on each (irreducible) component of $C$. An (irreducible) component of $C$ is called a ghost component
if the restriction $f$ to that component represents a trivial homology class.
Let $\nu$ be a section of the bundle $\Hom(\pi_1^*T\P^N,\pi_2^*TX)$ over
$\P^N\times X$ satisfying $J\circ \nu =-\nu\circ  J_{\P^N}$.
A stable map $f$ is $(J,\nu)$-holomorphic if
\begin{equation*}\label{hol-eq}
\frac12 (df +Jdfj)\, = \,(f,\phi)^*\nu
\end{equation*}
where  $j$ is the complex structure on $C$.
Denote by $\CM_{g,n}(X,A)$
the moduli space of $(J,\nu)$-holomorphic maps from nodal curves of (arithmetic) genus $g$ with $n$ marked points
that represent the class $A$ (we often omit $(J,\nu)$ in notation).
We also denote by
\begin{equation*}\label{moduli-ndg}
\CM_{\chi,n}^*(X,A)
\end{equation*}
the moduli space of $(J,\nu)$-holomorphic maps $f$ from possibly disconnected domains of Euler characteristic $\chi$
with {\em no degree zero connected components},
namely the restriction of $f$ to each ``connected component'' of its domain represents
a non-trivial homology class.

\medskip
For a finite set $A$, let $|A|$ denote the number of elements of $A$.

\begin{rem}\label{ghost}
A stable map $f$ in the moduli space $\CM_{\chi,n}^*(X,A)$ might have ghost components.
Let $C=C_1\cup C_2$ be the domain of $f$ such that $C_1$ is a connected curve that  is a union of some ghost components of $f$.
Then, the stability of $f$ implies that $2g(C_1) + \ell +n_1\geq 3$ where $\ell=|C_1\cap C_2|$ and
$n_1$ is the number of marked points on $C_1$.
\end{rem}

\medskip
The discussion below will be frequently used in  subsequent arguments.
Let $\phi_i$ be the first Chern class of line bundle over $\CM_{g,n}$ whose fiber over $(C,\{x_i\})$  is $T_{x_i}^*C$.
For a subset $I$ of $\{1,\cdots,n\}$,  let $\delta_I$ denote the Poincar\'{e} dual of the fundamental class of
the boundary stratum of $\CM_{g,n}$ that consists of  nodal curves $C_1\cup C_2$ where $C_1$ has genus zero,
$C_2$ has genus $g$ and the marked points on $C_1$ are precisely those labeled by $I$.

Consider the forgetful map
\begin{equation*}\label{forgetfulmap}
\pi_k : \CM_{g,n+k}\ \to\  \CM_{g,n}\ \ \ \ \
\end{equation*}
that forgets the last $k$ marked points.
For $1\leq i\leq n$, we have
\begin{equation}\label{AC1}
\pi_k^*\phi_i\,=\,\phi_i\,-\,\sum\delta_{\{i\}\cup I}
\end{equation}
where the sum is over all $I\subset \{n+1,\ldots,n+k\}$ with $I\ne \emptyset$ (cf. Lemma 3.1 of \cite{AC}).
The standard gluing map
\begin{equation}\label{gluingmap}
\eta : \CM_{g_1,n_1+1}\times \CM_{g_2,n_2+1}\ \to\ \CM_{g_1+g_2,n_1+n_2}
\end{equation}
is obtained by identifying  the last marked point of the first component with the first marked point on the second
component.
For our purpose, we extend this gluing map to the cases where $2g_2+n_2 <2$.
Denote  by $\CM_{g,n}$ the space of one point when $2g+n<3$ and note that
$$
\CM_{g,n_1+1}\times \CM_{0,1}\ \cong \ \CM_{g,n_1}
\ \ \ \ \ \ \ \mbox{and}\ \ \ \ \ \ \
\CM_{g,n_1+1}\times \CM_{0,2}\ \cong\  \CM_{g,n_1+1}.
$$
Let $\eta$ be the forgetful map that forgets the last marked point
when $2g_1+n_1\geq 3$ and $(g_2,n_2)=(0,0)$,  and let $\eta$ be the identity map when $(g_2,n_2)=(0,1)$.
The following fact then directly follows from (\ref{AC1}) and Lemma 3.3 of \cite{AC}.

\begin{lemma}[\cite{AC}]\label{L:AC}
For $1\leq i\leq n_1$ and $I\subset \{1,\cdots,n_1\}$ with $2\leq |I|<n_1$, we have
\begin{itemize}
\item[(a)]
if\  $2g_1+n_1\geq 3$ and $(g_2,n_2)=(0,0)$ then
\begin{equation*}\label{eta-R}
\eta^*\phi_i\, =\, (\phi_i - \delta_{\{i,n_1+1\}})\otimes 1
\ \ \ \ \ \ \ \mbox{and}\ \ \ \ \ \ \
\eta^*\delta_I\,=\,(\delta_I  + \delta_{I\cup \{n_1+1\}})\otimes 1
\end{equation*}
\item[(b)]
if either $(g_2,n_2)=(0,1)$ or  $2g_j+n_j\geq 2$ for $j=1,2$ then
\begin{equation*}\label{AC2}
\eta^*\phi_i\,=\,\phi_i\otimes 1
\ \ \ \ \ \ \ \mbox{and}\ \ \ \ \ \ \
\eta^*\delta_I\,=\,
\delta_I\otimes 1.
\end{equation*}
\end{itemize}
\end{lemma}

Lastly, we denote  the space of  curves with finitely many connected
components, Euler class $\chi$ and $n$ marked points by
$$
\TM_{\chi,n}
$$
(cf. page 57 of \cite{IP1}).
This space is a disjoint union of the products of the spaces $\CM_{g_j,n_j}$ with $\sum (2-2g_j)=\chi$ and $\sum n_j=n$
(including the unstable cases $2g_j+n_j<3$).
One can thus define $\phi_i$ classes and the boundary classes $\delta_I$ of $\TM_{\chi,n}$ in an obvious way.

\vskip 1 cm

\setcounter{equation}{0}
\section{Local GT Invariants}
\label{S:LGW}
\bigskip

In this section, we introduce local Gromov-Taubes invariants of spin curves and set up notation for them.
We will follow the definitions and notation in \cite{RT2}, \cite{LT} and \cite{IP1}.
Let $\pi:N\to D$ be a theta characteristic on a smooth curve $D$ of genus $h$.
The canonical bundle of the total space $N_D$ of $N$ is then isomorphic to $\pi^*N$,
so the tautological section of $\pi^*N$ gives a holomorphic 2-form $\a$
on $N_D$ whose zero set is the zero section $D\subset N_D$ (cf. \cite{LP1}).
The projectivization
$$
\P_h\,=\,\P(N\oplus \O_D)
$$
is a ruled surface over $D$. For small $\epsilon> 0$ fix an isomorphism (denoted by $\Psi$)
from a neighborhood of the zero section of $\P_h$
to the $3\epsilon$-neighborhood of $D$ in $N_D$ taking the zero section of $\P_h$ to $D$.
Choose a bump function
$\b$ that is 1 on the $\epsilon$-neighborhood of $D$ in $N_D$ and vanishes on the complement of $2\epsilon$-neighborhood
of $D$. The pull-back 2-form $\Psi^*(\b\a)$ is then a well-defined 2-form on $\P_h$.
Fix a fiber $V$ of the ruled surface $\P_h$ and for small $\delta>0$
choose a bump function $\b_{\sv}$ that is 1 on the
complement of the $2\delta$-neighborhood of $V$ and vanishes on $\delta$-neighborhood of $V$.
The 2-form
\begin{equation}\label{2-form}
\a_{\sv}\,=\,\b_{\sv} \Psi^*(\b\a)
\end{equation}
then induces, by (\ref{Ja}), an almost complex
structure $J_{\sv}$ on $\P_h$.
Let $U\subset \P_h$ be the preimage of $\epsilon$-neighborhood of $D$ in $N_D$ under the isomorphism $\Psi$.
We also denote by $D$ the zero section of $\P_h$ and
let $S$ be the section class represented by the zero section $D$ of $\P_h$.

\medskip

\begin{lemma}\label{ILL-F}
Let $U$ and $J_{\sv}$ be as above.
Then every $J_{\sv}$-holomorphic map from a connected domain into $\overline{U}$
that represents the class $dS$ with $d\ne 0$ is, in fact, holomorphic and its image
lies entirely in $D$, i.e.
$$
\CM_{\chi,n}^*(\overline{U},dS,J_{\sv})\,=\,\CM_{\chi,n}^*(D,d).
$$
\end{lemma}

\pf This proof is similar to that  of Lemma 3.2 in \cite{LP1}.
Use the same notation $\a$  for the real part of the  holomorphic 2-form $\a$
on $N_D$.
Let $f:(C,j)\to \overline{U}$ be a $J_{\sv}$-holomorphic map from a connected curve $C$ with complex structure $j$
that represents the class $dS$ where $d\ne 0$.  For each point $q\in C$, let $\{e_1,e_2=je_1\}$
be an orthonormal  basis of $T_qC$. Then
\begin{equation}\label{key-inequality}
 |\del f|^2 \,=\,   f^*|\b_{\sv}\Psi^*(\b\a)|^2|\bd f|^2 \,=\,
 f^*(\b_{\sv}\Psi^*(\b\a))(e_1,e_2) \,\leq\,  f^*(\Psi^*\a)(e_1,e_2)
\end{equation}
where the two equalities follow from Proposition 1.3 of \cite{L} and
the  inequality follows from the facts (i) $\Psi^*\b\equiv 1$ on $\overline{U}$ and
(ii) $0\leq \b_{\sv}\leq 1$.
Since $\a$ is a real part of holomorphic 2-form, integrating over the domain
shows $f$ is indeed holomorphic and  the image of $f$ lies in the zero set of
$\b_{\sv}\Psi^*(\beta\a)$ in $\overline{U}$.
Since $f$ represents the class $dS$, the image of $f$ must lie in the zero set of $\Psi^*(\beta\a)$ in $\overline{U}$
which is the zero section $D$. This completes the proof. \qed

\medskip
\begin{rem}\label{bump-ILL}
Let $\a$ be a holomorphic 2-form on $U$ and $g$ be any function on $U$ that satisfies $0\leq g\leq  1$.
Then, the inequality in (\ref{key-inequality}) shows that
for $\widetilde{\a}=g\a$ every $J_{\widetilde{\a}\,}$-holomorphic map into $U$ representing the class $dS$ ($d\ne 0$)
is holomorphic and has its image lying in the zero set of $\widetilde{\a}$.
\end{rem}

\bigskip

Let $d\ne 0$ and fix $(J,\nu)$ that is close to $(J_{\sv},0)$.
Lemma~\ref{ILL-F} and the Gromov Compactness Theorem imply that
the moduli space of $(J,\nu)$-holomorphic maps $\CM_{\chi,n}^*(U,dS)$ is compact.
The construction of Li and Tian \cite{LT} then defines the (virtual) fundamental class
\begin{equation}\label{LT}
\big[\,\CM_{\chi,n}^*(U,dS)\,\big]^{\vir}\ \in\ H_{*}\big(\,\Map_{\chi,n}(\P_h,dS)\,;\,\Q\,\big)
\end{equation}
in the homology of the space $\Map_{\chi,n}(\P_h,dS)$ of stable maps into $\P_h$
from nodal curves of Euler characteristic $\chi$ with $n$ marked points that
represent the homology class $dS$.

\begin{defn}\label{D:LGW}
The local (descendent) GT invariants of the spin curve $(D,N)$ of genus $h$ with parity $p\equiv h^0(N)$ (mod 2) are\,:
\begin{equation}\label{d-lgw}
GT_d^{loc,h,p}\big(\,\prod_{i=1}^n\tau_{k_i}(F^*)\,\big)\,=\,
\big[\,\CM_{\chi,n}^*(U,dS)\,\big]^{\vir}\, \cap\, \big(\,\prod_{i=1}^n \psi_i^{k_i}\cup ev_i^*(F^*)\,\big)
\end{equation}
where
$\psi_i$ be  the Euler class of the bundle over  $\Map_{\chi,n}(\P_h,dS)$  whose fiber
over $(f,C,\{x_i\})$ is $T^*_{x_i}C$,
$F^*$ is the Poincar\'{e} dual of the fiber class of $\P_h$,
$ev_i$ is the evaluation map at the $i$-th marked point,
and the Euler characteristic $\chi$ satisfies
\begin{equation*}\label{dim-d}
\sum k_i\,=\,d(1-h)\,-\,\tfrac12\,\chi.
\end{equation*}
We will often write respectively $+$ and $-$ for $p\equiv 0$ (mod 2) and for $p\equiv 1$ (mod 2).
\end{defn}

\begin{rem}\label{L:Inv-LGW}
Let $\b_t$ be a path  from $\b_0=\b_V$ to $\b_1\equiv 1$ with $0\leq \b_t\leq 1$ on $\P_h$ and let $J_t$
denote the almost complex structure induced from the 2-from $\b_t\Psi^*(\b\a)$ by (\ref{Ja}).
The proof of Lemma~\ref{ILL-F} shows for $d\ne 0$ and for all $t$
$$
\CM_{\chi,n}^*(\overline{U},dS,J_t) \,=\, \CM_{\chi,n}^*(D,d).
$$
In particular, this shows
$\CM_{\chi,n}^*(\overline{U},dS,J_t)$ is compact for all $t$.
It then follows from the standard cobordism argument (cf. Proposition 2.3 of \cite{LT})  that
the (virtual) fundamental class (\ref{LT}) is independent of the choice of $J_t$.
So, when $t=1$,  the isomorphism $\Psi$ as above gives
$$
\big[\,\CM_{\chi,n}^*(U,dS)\,\big]^{\vir}\,\cap\,\big(\,\prod_{i=1}^n \psi_i^{k_i}\cup ev_i^*(F^*)\,\big)
\,=\
\big[\,\CM_{\chi,n}^*(N_D,d[D])\,\big]^{\vir}\,\cap\,\big(\,\prod_{i=1}^n \psi_i^{k_i}\cup ev_i^*(\pi^* \gamma^*)\,\big)
$$
where $\gamma^*\in H^2(D)$ is Poincar\'{e} dual of the point class of $D$.
Thus, given $\chi,n$ and $d\ne 0$ the  invariant (\ref{d-lgw}) depends only on
the genus $h$ of $D$ and the parity  $p=h^0(N)$ (mod 2).
\end{rem}

\medskip

The stabilization and  evaluation at marked points defines a map
\begin{equation}\label{st-ev}
\ve\,=\,st\times ev\,:\,\Map_{\chi,n}(\P_h,dS)\ \to\ \TM_{\chi,n}\times (\P_h)^n.
\end{equation}
For the classes $\phi_i$ on $\TM_{\chi,n}$, we set
$$
GT_d^{loc,h,p}\big(\,\prod_{i=1}^n\phi_i^{k_i}(F^*)\,\big)\,=\,
\big[\,\CM_{\chi,n}^*(U,dS)\,\big]^{\vir}\, \cap\, \big(\,\prod_{i=1}^n st^*\phi_i^{k_i}\cup ev_i^*(F^*)\,\big).
$$
Now, suppose $(J,\nu)$ is generic (see page 10 of \cite{RT2}).  Then,
the image of  $\CM_{\chi,n}^*(U,dS)$ under the map (\ref{st-ev}) defines a homology class
\begin{equation}\label{RTclass}
\big[\,\CM_{\chi,n}^*(U,dS)\,\big]\ \in\ H_{*}\big(\,\TM_{\chi,n}\times (\P_h)^n\,;\,\Q\,\big)
\end{equation}
satisfying
$\ve_* \big[\,\CM_{\chi,n}^*(U,dS)\,\big]^{\vir}= \big[\,\CM_{\chi,n}^*(U,dS)\,\big]$
(cf. Remark 10.2 of \cite{LP2}). So, we have
\begin{equation*}\label{LT=RT}
GT_d^{loc,h,p}\big(\,\prod_{i=1}^n\phi_i^{k_i}(F^*)\,\big)\,=\,
\big[\,\CM_{\chi,n}^*(U,dS)\,\big]\, \cap\, \prod_{i=1}^n \phi_i^{k_i}\otimes (F^*)^{\otimes n}.
\end{equation*}

\medskip
\begin{rem}\label{h>0}
If the spin curve $(D,N)$ has genus $h>0$, then
there are no non-constant holomorphic maps
from genus zero curves to $D$ and hence, by the Gromov Compactness Theorem and Lemma~\ref{ILL-F}, for every map $f$ in the moduli space
$\CM_{\chi,n}^*(U,dS)$ every genus zero (irreducible) component is ghost component. It thus follows from the stability and the relation between
$\psi_i$ class and $st^*\phi_i$ class (cf. \cite{KM} page 388) that
\begin{equation}\label{t=p}
GT_d^{loc,h,p}\big(\,\prod_{i=1}^n\tau_{k_i}(F^*)\,\big)\,=\,
GT_d^{loc,h,p}\big(\,\prod_{i=1}^n\phi_i^{k_i}(F^*)\,\big).
\end{equation}
\end{rem}

\medskip
We end this section with  dimension zero local GT invariants for $d=1$ and 2.

\medskip

\begin{lemma}\label{com:deg=1,2}\ \ \ \ \
$GT_1^{loc,h,p}\,=\,(-1)^p$\ \ \ \ and\ \ \ \
$GT_2^{loc,h,p}\,=\,(-1)^p 2^{h-1}$.
\end{lemma}

\pf
The dimension zero  local GW invariants $GW_d^{loc,h,p}$   and
the dimension zero local GT invariants $GT_d^{loc,h,p}$ are related as follows\,:
$$
1\,+\,\sum_{d>0}GT_d^{loc,h,p}\,t^d\,=\,\exp{\Big(\sum_{d>0} GW_d^{loc,h,p}\,t^d\Big)}
$$
(cf. Section 2 of \cite{IP2}). The lemma thus follows from the fact
$$
GW_1^{loc,h,p}\,=\,(-1)^p\ \ \ \ \ \ \
\mbox{and}\ \ \ \ \ \ \
GW_2^{loc,h,p}\,=\,\frac12\big[ (-1)^p\,2^h-1\big]
$$
(see Section 10 of \cite{LP1}). \qed

\vskip 1cm
\setcounter{equation}{0}
\section{Relative Local Invariants}
\label{S:RLGW}
\bigskip

In \cite{IP1}, GW invariants were generalized to relative GW invariants relative to codimension two
symplectic submanifold. Following \cite{IP1}, we can define relative local invariants.
A  $(J,\nu)$-holomorphic map $f$ is called {\em $V$-regular with a contact vector
$s=(s_1,\ldots,s_\ell)$ } if $f^{-1}(V)$ consists of
the last $\ell$ ordered marked points $x_{n+1},\ldots,x_{n+\ell}$ such that the image of $f$ has
the contact order $s_k$ at $x_{n+k}$. Denote by
\begin{equation*}\label{rel-moduli}
\M_{\chi,n,s}^{V,*}(U,dS)
\end{equation*}
the moduli space of $V$-regular  $(J,\nu)$-holomorphic maps $f$
into $U$ with contact vector $s$ where the superscript $*$ also means $f$ has no degree zero connected components.
For a contact vector $s=(s_1,\ldots,s_\ell)$, we write
\begin{equation*}\label{deg-length}
\deg(s)=\sum_{k=1}^\ell s_k,\ \ \ \ \ \ \ \ \ \ \
\ell(s)=\ell,\ \ \ \ \ \ \ \ \ \ \
|s|=\prod_{k=1}^\ell s_k
\end{equation*}
and,  noting there are no rim tori  since $H_1(V)=0$ (cf. Remark 5.3 of \cite{IP1}), we set
$$
V_s\ =\ \{\ \big(\,(v_1,s_1),\ldots,(v_\ell,s_\ell)\,\big)\ \,|\,\ v_k\in V\ \}.
$$
The moduli space of $V$-regular maps $\M_{\chi,n,s}^{V,*}(U,dS)$ also has an associated map
\begin{equation}\label{st-ev-h}
\ve_s=st\times ev\times h_s\,:\,\M_{\chi,n}^{V,*}(U,dS)\,\to\, \TM_{\chi,n+\ell(s)}\times (\P_h)^n\times  V_s
\end{equation}
where $ev$ is the evaluation map at the first $n$ marked points and $h_s$ is given by
\begin{equation}\label{ev-h}
h_s(f,x_1,\cdots,x_{n+\ell})\  = \ \big(\,(f(x_{n+1}),s_1),\cdots,(f(x_{n+\ell}),s_\ell)\,\big).
\end{equation}

Observe that the (holomorphic) fiber $V$ of $\P_h$ is $J_{\sv}$-holomorphic since the 2-form $\a_{\sv}$ in (\ref{2-form}) vanishes near $V$.
The pair $(J_{\sv},0)$ is thus $V$-compatible in the sense of Definition 3.2 of \cite{IP1}.
Now, choose a  generic $V$-compatible  $(J,\nu)$ that is sufficiently close to $(J_{\sv},0)$.
Lemma~\ref{ILL-F}, the Gromov Compactness Theorem and the relative GW theory of \cite{IP1} then imply that
the image of the moduli space  $\M_{\chi,n,s}^{V,*}(U,dS)$
under the map (\ref{st-ev-h})
defines a homology class
\begin{equation}\label{IP}
\big[\,\M^{V,*}_{\chi,n,s}(U,dS)\,\big]\ \in\
 H_{*}\big(\,\TM_{\chi,n+\ell(s)}\times (\P_h)^n\times V_s\,;\,\Q\,\big)
\end{equation}

Let  $\{\b_j\}$ be a basis of $H_*(V)$. Then a basis of $H^*(V_s)$ is given by elements of the form
$$
C_{J,s}\,=\,C_{\beta_{j_1},s_1}\otimes\cdots\otimes C_{\beta_{j_\ell},s_\ell}.
$$
Lemma~\ref{ILL-F} and the Gromov Compactness Theorem imply that
\begin{equation}\label{pt-constraint}
\big[\,\M_{\chi,n,s}^{V,*}(U,dS)\,\big]\,\cap\,C_{J,s}\,=\,0
\end{equation}
unless all $\b_{j_k}$ are the fundamental class $[V]$ of $V$.
Since $H_1(V)=0$, we can forget the ordering of the contact constraints $C_{J,s}$
by simply writing
\begin{equation}\label{contact-no-order}
\prod_{j,b}\,(C_{\b_j,b})^{m_{j,b}}\, =\, C_{\beta_{j_1} ,s_1}\cdots C_{\beta_{j_\ell} ,s_\ell}
\end{equation}
with the relation $C_{\b_j,b}\cdot C_{\b_i,a}= C_{\b_i,a}\cdot C_{\b_j,b}$ where
$m=(m_{j,b})$ is a sequence of nonnegative integers determined by  (\ref{contact-no-order}).
If all $\b_{j_k}=\b_j$ for some $j$ then the sequence $m$ can be considered as a partition
of the integer $d$.
The (unordered) contact constraint (\ref{contact-no-order}) is then a pair of the Poincar\'{e} dual of $\b_j$ and
the partition $m$ of $d$, i.e. $m=(m_1,\cdots,m_\ell)$ with $m_1\leq m_2\leq\cdots\leq m_\ell$ and
$\sum m_j=d$. Write $m=(1^d)$ if all $m_j=1$. In that case,
$$
|m|\,=\,|(1^d)|\,=\,1\ \ \ \ \ \mbox{and}\ \ \ \ \ m!\,=\,(1^d)!=d!.
$$
When all  $\b_{j_k}=\b_j$ for some $j$, we write  $C_{J,s}$  simply as $C_{\b_j^{\ell}}$.

\medskip
\begin{rem}\label{partition}
Given a partition $m=(m_1,\cdots,m_\ell)$ of $d$, there are $\ell!/m!$ ordered sequences $s$ with
$\tau(s)=m$ for some permutation $\tau$ in the symmetric group $S_\ell$.
\end{rem}

\medskip

\begin{defn}
For a partition $m$ of $d$ with $m=\tau(s)$ for some permutation $\tau$ in $S_{\ell(s)}$,
we set
\begin{equation*}\label{r-lgw}
GT_{m}^{loc,h,p}\big(\,\prod_{i=1}^n\phi_i^{k_i}(F^*)\,\big)
\,=\,
\big[\,\M_{\chi,n,s}^{V,*}(U,dS)\,\big]\,\cap\,
\prod_{i=1}^n \phi_i^{k_i}\otimes (F^*)^n \otimes C_{[V]^{\ell(m)}}
\end{equation*}
where the Euler characteristic $\chi$ is given by
\begin{equation*}\label{dim-r}
\sum k_i\,=\,d(1-h)\,-\,\tfrac12\,\chi\,+\,(\ell(m)-d).
\end{equation*}
\end{defn}

Choose two distinct fibers $V_1$ and $V_2$ of $\P_h$ and let $V=V_1\sqcup V_2$.
We can also define an almost complex structure $J_{\sv}$
that equals to the K\"{a}hler structure of $\P_h$ near  $V$ and
satisfies Lemma~\ref{ILL-F}. Thus, we can  define relative local invariants
\begin{equation}\label{mul-rel}
GT_{m^1,m^2}^{loc,h,p}\big(\,\prod_{i=1}^n\phi_i^{k_i}(F^*)\,\big)
\end{equation}
relative to $V$ (with contact vectors $m^i$ with $V_i$) for  the class $dS$ with Euler characteristic $\chi$ where
the Euler characteristic $\chi$ satisfies
$$
\sum k_i\,=\,d(1-h)\,-\,\tfrac12\,\chi\,+\,\sum (\ell(m^i)-d).
$$

\begin{rem}\label{h=0}
The only genus zero spin curve is the even spin curve $(\P^1,\O(-1))$. In this case,
 $S^2=-1$ and hence for $d\ne 0$ we have
\begin{equation*}\label{e:h=0}
\CM_{\chi,n}^*(U,dS,J_{\sv})\,=\,\CM_{\chi,n}^*(\P_0,dS,J_{\sv}).
\end{equation*}
This shows that degree $d$ local invariants of spin curve $(\P^1,\O(-1))$ are the same as the GW invariants of
$\P_0$ for the class $dS$. It also shows that relative local invariants are  equal to
the relative GW invariants of $(\P_0,V)$.
\end{rem}

For simplicity, we set
$$
\CM\,=\,\CM_{\chi,n}^*(U,dS)\ \ \ \ \
\mbox{and}\ \ \ \ \
\M^{V}\,=\,\M^{V_1,V_2,*}_{\chi,n,m^1\!,m^2}(U,dS).
$$
Noting the homology class (\ref{RTclass}) defines a  map
$H^{*}((\P_h)^n)\to H_{*}(\TM_{\chi,n})$, we set
\begin{equation*}
GT_{d,\chi}^{loc,h,p}\big(\,(F^*)^n\,\big)\, =\,
\big[\,\CM\,\big]\,\cap\,(F^*)^{\otimes n}\,\in\,H_{2q}(\TM_{\chi,n})
\end{equation*}
where $q=d(1-h)-\frac12\chi$. Similarly, we also set
$$
GT_{m^1,m^2,\chi}^{loc,h,p}\big(\,(F^*)^n\,\big)\,=\,
[\,\M^{V}\,]\,\cap\,(F^*)^{n}\,\bigotimes_{i=1}^2\, C_{[V_i]^{\ell(m^i)}}
\,\in\,H_{2r}(\TM_{\chi,n+\sum \ell(m^i)})
$$
where $r=d(1-h)-\tfrac12\chi +\sum (\ell(m^i)-d)$.

\medskip

\begin{rem}\label{cut-down-class}
Let $B$ be a geometric representative of the $n$ product of fiber classes $F^{\otimes n}$ of $(\P_h)^n$
in general position with respect to the
evaluation map at marked points. Then the images of the cut-down moduli spaces $\CM \cap B$ and
$\M^{V}\cap B$ under the stabilization map respectively define classes satisfying
$$
\big[\,st\big(\,\CM\cap B\,\big)\,\big]\,=\,
GT_{d,\chi}^{loc,h,p}\big(\,(F^*)^n\,\big)
\ \
\mbox{and}\ \
\big[\,st\big(\,\M^{V}\cap B\,\big)\,\big]\,=\,
GT_{m^1,m^2,\chi}^{loc,h,p}\big(\,(F^*)^n\,\big).
$$
\end{rem}

\medskip
\medskip
For the ruled surface $\F_0=\P^1\times\P^1$, we use the same notations
$F$ and $S$ for the fiber class and the section class, respectively.
To save notation, we also use the same notation $V$  for the union of $2$ distinct fibers
$V_1$ and $V_2$ of the ruled surface $\F_0$. For partitions $m^i$ of $d$, denote by
\begin{equation}\label{rel-F0}
GT^{\F_0}_{m^1,m^2}\big(\,\prod_{i=1}^n \phi_i^{k_i}(F^*)\,\big)
\end{equation}
the relative GT invariants of  $(\F_0,V)$
for the class $dS$ with Euler characteristic $\chi$ satisfying
$$
\sum k_i\,=\,d\,-\,\tfrac12\,\chi\,+\,(\,\ell(m^2)-d\,)
$$ where the contact constraint with $V_1$ and $V_2$ are respectively
$C_{pt^{\ell(m^1)}}$ and  $C_{[V_2]^{\ell(m^2)}}$.
It follows directly  from Lemma 14.6 of \cite{IP2} that
\begin{equation}\label{com-F0}
GT^{\F_0}_{(1^d)}\,=\,1
\ \ \ \ \ \ \ \mbox{and}\ \ \ \ \ \ \
GT^{\F_0}_{(1^d),(1^d)}\,=\,d!
\end{equation}
For the class  $\big[\,\M^{V}_{\chi,n}(\F_0,dS)\,\big]$ that defines
the relative invariants (\ref{rel-F0}) of $(\F_0,V)$ we also set
$$
GT^{\F_0}_{m^1,m^2,\chi}\big(\,(F^*)^n\,\big)\,=\,
\big[\,\M^{V}_{\chi,n}(\F_0,dS)\,\big]
\,\cap\,(F^*)^{n}\otimes C_{pt^{\ell(m^1)}}\, \otimes\, C_{[V_2]^{\ell(m^2)}}.
$$
This is a homology class in $H_{2t}(\TM_{\chi,n+\sum \ell(m^i)})$ where $t=d-\tfrac12\chi +\sum (\ell(m^2)-d)$.

\vskip 1cm

\setcounter{equation}{0}
\section{Blow-Up and Sum Formula}
\label{S:Blow-Up}
\bigskip

The aim of this section is to prove Theorem~A in the Introduction.
We will apply the limiting and smoothing arguments of \cite{IP2}
to  our local invariants.
The proof consists of three steps.

\bigskip
\non
{\bf Step 1 :}
Fix a fiber $V_0$ of $\P_h$ and consider a degeneration
\begin{equation*}
\lam \,:\, Z\,\overset{\sigma}{\longrightarrow}\,
 \P_h\times \cx \,\longrightarrow \,\cx
\end{equation*}
where $\sigma:Z\to \P_h\times \cx$ is the blow-up of $\P_h\times \cx$ along $V_0\times \{0\}$ and the second map
is projection onto the second factor.
The central fiber $Z_0$ is the singular surface $\P_h \cup_{V_0} \F_0$ and general fibers
$Z_\lam$ ($\lam\ne 0$) are isomorphic to $\P_h$ that is the symplectic fiber sum $\P_h\#_{V_0}\F_0$.
For $\lam\ne 0$, let $D_\lam$ denote the zero section of $Z_\lam$, i.e., $D_\lam=\sigma^{-1}(D\times \{\lam\})$.

Choose a fiber $V_1\ne V_0$ of $\P_h\subset Z_0$  and a fiber $V_2\ne V_0$ of $\F_0\subset Z_0$
and set
$$
V\,=\,V_1\sqcup V_0\sqcup V_2.
$$
One can choose a (smooth) family $\tilde{V}_\lambda$ of disjoint union of two fibers of $Z_\lam\simeq \P_h$
with $\tilde{V}_0=V_1\sqcup V_2$.
Denote by $\J(Z)$ the space of all $(J,\nu)$ on $Z$ satisfying (i)
each $Z_\lam$ is $J$-invariant and (ii)
the restriction of $(J,\nu)$ to $Z_\lam$ ($\lam\ne 0$) is $\tilde{V}_\lam$-compatible
and to $Z_0$ is $V$-compatible
(cf. Lemma 2.3 of \cite{IP2}).
We will  use the same notation $(J,\nu)$ for its restriction to each $Z_\lam$.

Fix a small $\delta>0$ and define a $\delta$-neck $Z(\delta)$ as a (normal) $\delta$-neighborhood of $V$ in $Z$.
The energy of a map $f$ (more precisely of $(f,\phi)$ as in (\ref{energy})) into $Z$ in the $\delta$-neck is
\begin{equation*}\label{energy-neck}
E^\delta(f) = \frac12\int |df|^2 +|d\phi|^2
\end{equation*}
where the integral is over $f^{-1}(Z(\delta))$.
By Lemma\,1.5 of \cite{IP1} there is a constant $c_{\sv}$ depending
only on the restriction of $(J,\nu)$ to $V\subset Z$ such that every component of every  $(J,\nu)$-holomorphic map
into $V$ has energy greater that $c_{\sv}$.
A $(J,\nu)$-holomorphic map $f$ into $Z$ is {\em $\delta$-flat} if the energy
in the $\delta$-neck $E^\delta(f)$ is at most $c_{\sv}/2$. Note that a $\delta$-flat map into $Z_0$ has
no components mapped entirely into $V$.

Once and for all, fix $\chi, n, d\ne 0$ and
for each $\lam\ne 0$  we set
$$\M(Z_\lam)=\M_{\chi,n,(1^d),(1^d)}^{\tilde{V}_\lam}(Z_\lam,dS).$$
Denote the set of $\delta$-flat maps in $\M(Z_\lam)$  by $\M^\delta(Z_\lam)$ and write
\begin{equation}\label{limit-set}
\lim_{\lam\to 0}\,\M^\delta(Z_\lam)
\end{equation}
for the set of limits of sequences of $\delta$-flat maps in $\M^\delta(Z_\lam)$ as $\lam\to 0$.
Since $\delta$-flatness is a closed condition,
each map $f$ in the limit set (\ref{limit-set}) is also $\delta$-flat and hence the domain of $f$
has no components mapped entirely into $V$. Consequently, we have
\begin{itemize}
\item[$(a)$] $f$ splits  as $f=(f_1,f_2)$
where $f_1$  and $f_2$ are respectively  $(V_1\sqcup V_0)$-regular map into $\P_h$ and
$(V_0\sqcup V_2)$-regular map into $\F_0$ and each $f_i$ has contact vector $(1^d)$ with $V_i$ for $i=1,2$,
\item[$(b)$] $f^{-1}(V_0)$ consists of nodes $\{p_1,\cdots,p_\ell\}$ of the domain
such that each $p_i$ has a well-defined multiplicity $s_i$ equal to the order of contact of the image of $f_1$ (or $f_2$)
with $V_0$ at $p_i$
\end{itemize}
(see Section 3 of \cite{IP2}).
Renumbering the nodes $\{p_1,\cdots,p_\ell\}$  gives $\ell!$ ordered sequences $s=(s_1,s_2,\cdots,s_\ell)$.
On the other hand, since for small $|\lam|$ the $\delta$-flat maps in $\M^\delta(Z_\lambda)$ are $C^0$-close to
$\delta$-flat maps
in the limit set (\ref{limit-set}), to each map $f$ in $\M^\delta(Z_\lambda)$ one can assign ordered sequences
$s$.  Denote by $\M_s^\delta(Z_\lam)$ the set of all such pair $(f,s)$ labeled by an ordered sequence $s$.
Then, there are actions of symmetric groups $S_{\ell}$ such that
\begin{equation}\label{action}
 \bigsqcup_{\ell}\,\Big(\bigsqcup_{\ell(s)=\ell} \,\,\M_s^\delta(Z_\lam)\Big)\,\Big/ S_\ell
\ =\ \M^\delta(Z_\lam).
\end{equation}

For $\P_h,\F_0\subset Z_0$ and each ordered sequence $s$ with $\deg(s)=d$, there is an evaluation map
\begin{equation}\label{evaluation}
ev_s\,:\,{\textstyle \bigcup}\,
\Big(\,\M_{\chi_1,n_1,(1^d),s}^{V_1,V_0}(\P_h,dS-kF)
\times
\M_{\chi_2,n_2,s,(1^d)}^{V_0,V_2}(\F_0,dS+kF)
\,\Big)
\longrightarrow\ V^{\ell(s)}_0\times V^{\ell(s)}_0
\end{equation}
that records the intersection points with the fiber $V_0$
where the union is over all $0\leq k\leq h$, $n_1+n_2=n$ and
$\chi=\chi_1+\chi_2-2\ell(s)$. Let $\triangle_s$  be the diagonal of $V_0^{\ell(s)}\times V_0^{\ell(s)}$ and
denote by
$$
\K_s^\delta\ \subset\ ev_s^{-1}(\triangle_s)
$$
the set of $\delta$-flat maps in $ev_s^{-1}(\triangle_s)$.
Since each map in $ev_s^{-1}(\triangle_s)$
can be considered as  a  pair of a map $f$ into $Z_0=\P_h\cup_{V_0} \F_0$
satisfying $(a)$ and $(b)$ with an ordered sequence $s$, we have
\begin{equation}\label{limit}
\lim_{\lam\to 0}\ \bigsqcup_s\ \M_s^\delta(Z_\lam)\ \subset\ \bigsqcup_s\,\K_s^\delta.
\end{equation}

Conversely,  each  map $f=(f_1,f_2)\in \K_s^\delta$ can be smoothed to
produce  maps in $\M_s^\delta(Z_\lam)$ for small $|\lam|$.
Let $C_1$ and $C_2$ be  the domains of $f_1$ and $f_2$ respectively.
Identifying  the $\ell(s)$ contact points with $V_0$ of $C_1$ with
the $\ell(s)$ contact points with $V_0$  of $C_2$
determines a gluing map
$$
\TM_{\chi_1,n_1+\ell(s)+d}\times \TM_{\chi_2,n_2+\ell(s)+d}\ \longrightarrow\
\TM_{\chi,n+2d}
$$
where $\chi_i=\chi(C_i)$ for $i=1,2$. For each $\ell(s)$,
taking the union over all $\chi_1,\chi_2,n_1$ and $n_2$ defines a gluing map
\begin{equation*}\label{attaching}
\xi_{\ell(s)} \ :\  \bigsqcup\, \TM_{\chi_1,n_1+\ell(s)+d}\times \TM_{\chi_2,n_2+\ell(s)+d} \ \longrightarrow\
\TM_{\chi,n+2d}
\end{equation*}
Theorem 10.1 of \cite{IP2} then gives\,:

\begin{theorem}[\cite{IP2}]\label{T:smoothing}
For generic $(J,\nu)\in\J(Z)$ and for small $|\lam|$, there is an  $|s|$-fold covering
$$
\pi_{s,\lam}:\K_s^\delta(Z_\lam)\to \K^\delta_s
$$ with a commutative diagram (up to homotopy)\,:
\begin{equation}\label{smoothing}
\xymatrix{
 \ \ \ \ \ \ \ \ \ \ \ \ \ \ \bigsqcup_s\, \K_s^\delta(Z_\lam)\ \ \ \ \ \ \ \ \ \ \ \ \ \ \ \
 \ar[d]^{st_\lambda} \ar[rr]^{\ \ \ \ \ \ \Phi_\lam}
 &&  \ \ \ \ \ \ \ \ \bigsqcup_s\, \M_s^\delta(Z_\lam) \ar[d]^{st} \ \ \ \ \ \\
 \ \ \ \ \ \bigsqcup_s\,\TM_{\chi_1,n_1+\ell(s)+d}\times \TM_{\chi_2,n_2+\ell(s)+d} \ \ \ \ \ \ \
 \ar[rr]^{\ \ \ \ \ \ \xi}
 && \ \ \ \ \ \ \ \ \ \ \ \TM_{\chi,n+2d} \ \ \ \ \ \ \ \ \ \
}
\end{equation}
where the top arrow is an embedding, $st_\lambda=st\circ \pi_{s,\lam}$ and
$\xi$ in the bottom arrow is given by the gluing maps $\xi_{\ell(s)}$.
The construction of the smoothing map $\Phi_\lam$ also shows that
\begin{equation}\label{lim-Phi}
\lim_{\lam\to 0}\,\Phi_{\lam}(\tilde{f}_\lam)\ =\ f\ \ \ \ \ \ \
\mbox{where}\ \ \ \ \ \ \pi_{s,\lam}(\tilde{f}_\lam)=f.
\end{equation}
\end{theorem}

\bigskip

\non
{\bf Step 2 :}
Let $U$ be an open neighborhood of the zero section $D$ of $\P_h$
and $\a_0=\a_{\sv_0}$ be a 2-form
as in Section~\ref{S:LGW} (see the paragraph above Lemma~\ref{ILL-F}).
Regard the 2-form $\a_0$ as a 2-form on $\P_h\times \cx$ in an obvious way and for small $\epsilon>0$ choose
a bump function $\b$ that is 1 on the complement of $2\epsilon$-neighborhood of $V_1$ in $Z$ and
vanishes on $\epsilon$-neighborhood of $V_1$ in $Z$.
The 2-form $\b(p^*\a_0)$ then defines,
again by the formula (\ref{Ja}), an almost complex structure $J_{\sv}$ on $Z$.
The pair $(J_{\sv},0)\in \J(Z)$ and
the restriction of $J_{\sv}$ to $\F_0$
is the product complex structure of $\F_0$
since $p^*\a_0$ vanishes on some neighborhood of $\F_0$ in $Z$.
Moreover, by Remark~\ref{bump-ILL},  we have
\begin{align}\label{ILL-last}
\CM_{\chi,n+2d}^*(U_\lam,dS,J_{\sv})\,&=\,\CM^*_{\chi,n+2d}(D_\lam,d)_{\displaystyle\phantom{\sum}}
\notag \\
\CM_{\chi,n+2d}^*(U_0\cap \P_h,dS,J_{\sv})\,&=\,\CM^*_{\chi,n+2d}(D,d)
\end{align}
where $U_\lam =p^{-1}(U\times \{\lam\})$.
Fix an open neighborhood $W$ of $D$ in $\P_h$ satisfying $\overline{W}\subset U$ and for each $\lam$  set
$$
W_\lam = p^{-1}(W\times \{\lam\})\ \ \ \ \ \
$$
The following fact is our  key observation for the proof of Theorem~A.

\begin{lemma}\label{key}
For $(J,\nu)\in\J(Z)$ sufficiently close to $(J_{\sv},0)$ and for small $|\lam|>0$, we have
$$
\CM_{\chi,n+2d}^*(U_\lam,dS)\,\setminus\,\CM_{\chi,n+2d}^*(W_\lam,dS)\,=\,\emptyset.
$$
\end{lemma}

\pf Suppose not. Then there exists
a sequence of  $(J_k,\nu_k)$-holomorphic maps $f_k$ into $U_{\lam_k}$ with
$\im(f_k)\cap (U_{\lam_k}\setminus W_{\lam_k})\ne \emptyset$ and with
no degree zero connected components
where $\lam_k$ converges to $0$ and $(J_k,\nu_k)$
converges to $(J_{\sv},0)$ as $k\to \infty$.
The Gromov Compactness Theorem then implies that after passing to subsequences, $f_k$ converges to
a $J_{\sv}$-holomorphic map $f$ into $Z_0$ such that
(i) $\im(f)\subset \overline{U}\!_0$ and (ii) $\im(f)\cap (\overline{U}\!_0\setminus W_0)\ne\emptyset$.
Since the limit map $f$ also has no degree zero connected components, (i) implies
$f$ can be split as $f=(f_1,f_2)$ where
$f_1$ and $f_2$ map into $\P_h$ and $\F_0$  respectively such that
$$
\im(f_1)\,\cap\, V_0 \,=\,\im(f_2)\,\cap\, V_0.
$$
It follows from (\ref{ILL-last}) that $\im(f_1)\subset D$ and hence $\im(f_2)\cap V_0=D\cap V_0$.
Then, since the restriction of $J_{\sv}$ on $\F_0$
is the product complex structure, $\im(f_2)$ lies in the  section of
$\F_0$ passing through the intersection point $D\cap V_0$.  We have $\im(f)\subset W_0$
which contradicts (ii).
 \qed

\bigskip
Fix $(J,\nu)\in \J(Z)$ sufficiently close to $(J_{\sv},0)$ and for small $|\lam|$ set
$$
\M^{\delta,*}_s(U_\lam) \,=\,
\{(f,s)\in \M_s^\delta(Z_\lam)\,|\,f\in \M^{\tilde{V}_\lam,*}_{\chi,n,(1^d),(1^d)}(U_\lam,dS)\,\}.
$$
Consider the restriction of the evaluation map (\ref{evaluation})\,:
\begin{equation*}
ev_{s,U_0}^*\,:\,{\textstyle \bigcup}\,
\Big(\,\M_{\chi_1,,n_1,(1^d),s}^{V_1, V_0,*}(\P_h\cap U_0,dS)
\times
\M_{\chi_2,n_2,s,(1^d)}^{V_0, V_1,*}(\F_0\cap U_0,dS)
\,\Big)
\longrightarrow\ V^{\ell(s)}_0\times V^{\ell(s)}_0
\end{equation*}
where the union is over all $n_1+n_2=n$ and
$\chi=\chi_1+\chi_2-2\ell(s)$.

\begin{rem}
Let $q$ be the intersection point of $D$ and $V$ in
$Z_0=\P_h \cup_{V_0} \F_0$.
Then there is a unique section $E_q$ of $\F_0$ that lies in $\F_0\cap U_0$ and intersects with $V_0$ at the point $q$.
Choose $\ell(s)$ points $\{q_j\}$ in $V_0\cap (\F_0\cap U_0)$ that are sufficiently close to $q$. Denote by
$$
\M_{\chi_2,n_2,s,(1^d)}^{V_0, V_2,*}(\F_0,dS)\,\cap\, \{q_j\}
$$
the cut-down moduli space of $(V_0\sqcup V_2)$-regular $(J,\nu)$-holomorphic maps $f$  with $\ell(s)$ contact points
$\{x_{n+j}\}$ (with $V_0$) satisfying $f(x_{n+j})=q_j$.
Since every holomorphic map representing the class $dS$ and passing through the point $q$
has its image in $E_q$, by the Gromov Compactness Theorem,   we have
\begin{equation}\label{ILL-ruled-0}
\M_{\chi_2,n_2,s,(1^d)}^{V_0, V_2,*}(\F_0\cap U_0,dS)\,\cap\, \{q_j\}\,=\,
\M_{\chi_2,n_2,s,(1^d)}^{V_0, V_2,*}(\F_0,dS)\,\cap\, \{q_j\}.
\end{equation}
This shows that local invariants of $\F_0$ counting maps into $U_0$ with point constraints equal to
the standard invariants of $\F_0$ with points constraints.
\end{rem}

\medskip

Let
${\displaystyle
\K_{s,U_0}^{\delta,*}\, =\,  \K_s^\delta \ \cap\ (ev_{s,U_0}^*)^{-1}(\triangle_s)
}$.
Lemma~\ref{key} and (\ref{limit}) imply
\begin{equation}\label{limit-loc}
\begin{array}{c}
{\displaystyle
\lim_{\lam\to 0}\ \bigsqcup_s\,\M_s^{\delta,*}(U_\lam)\, \subset
\, \bigsqcup_{\phantom{A}s\phantom{A}}\,\K_{s,U_0}^{\delta,*}
}
\\
{\displaystyle
\lim_{\lam\to 0}\ \bigsqcup_s^{\phantom{s}}\,
\Big(\,\M_s^{\delta}(Z_\lam)\setminus \M_{s}^{\delta,*}(U_\lam)\,\Big)\
\  {\textstyle \bigcap}\ \ \bigsqcup_s\,\K_{s,U_0}^{\delta,*} \, =\, \emptyset.
}
\end{array}
\end{equation}
Consequently, by (\ref{smoothing}), (\ref{lim-Phi}) and (\ref{limit-loc}),
for the restriction  $\Phi_\lam^{loc}$  of the smoothing map $\Phi_\lam$ to
$$
\K_{s,U_0}^{\delta,*}(U_\lam)\,=\,\pi_{s,\lam}^{-1}\,(\K_{s,U_0}^{\delta,*})
$$
we have a commutative diagram (up to homotopy):
\begin{equation}\label{smoothing-local}
\xymatrix{
 \ \ \ \ \ \ \ \ \ \ \ \ \ \ \  \bigsqcup_s\, \K_{s,U_0}^{\delta,*}(U_\lam)\ \ \ \ \ \ \ \ \ \ \ \ \ \ \ \ \
 \ar[d]^{st_\lambda} \ar[rr]^{\ \ \ \ \ \ \ \ \Phi_\lam^{loc}}
 &&  \ \ \ \ \ \ \ \ \ \bigsqcup_s\, \M_s^{\delta,*}(U_\lam) \ar[d]^{st} \ \ \ \ \ \ \\
 \ \ \ \ \  \bigsqcup_s\,\TM_{\chi_1,n_1+\ell(s)+d}\times \TM_{\chi_2,n_2+\ell(s)+d} \ \ \ \ \ \ \
 \ar[rr]^{\ \ \ \ \ \ \ \ \xi }
 && \ \ \ \ \ \ \ \ \ \ \ \ \ \TM_{\chi,n+2d} \ \ \ \ \ \ \ \ \ \ \ \
}
\end{equation}

\bigskip
\non
{\bf Step 3 :}
The commutative diagram (\ref{smoothing-local})  leads to a sum formula for local invariants for the sum
$(\P_h,V_1\sqcup V_2)$ of
$(\P_h,V_1\sqcup V_0)$ and $(\F_0,V_0 \sqcup V_2)$ along $V_0$. We first assume that all maps in
\begin{equation}\label{msp-last}
\M^{\tilde{V}_\lam,*}_{\chi,n,(1^d),(1^d)}(U_\lam,dS)
\end{equation}
are $\delta$-flat when $|\lam|$ is small.
For fixed $n_1+n_2=n$, one can choose a continuous family of geometric representatives $B_\lam$ disjoint with $\tilde{V}_\lam$ satisfying\,:
\begin{itemize}
\item each $B_\lam$ ($\lam\ne 0$) is a geometric representative of the
$n$ product of fiber classes
$ F^{\otimes n}=F\otimes \cdots \otimes F$ of  $(Z_\lam)^n \simeq (\P_h)^n$,
\item $B_0=B_{\P_h}\sqcup B_{\F_0}$ where $B_{\P_h}$ and $B_{\F_0}$ are geometric representatives of
the classes $F^{\otimes n_1}$ of $(\P_h)^{n_1}$ and $F^{\otimes n_2}$  of $(\F_0)^{n_2}$ respectively.
\end{itemize}
It now follows from the diagram (\ref{smoothing-local}) that
\begin{align}\label{Thm10.4}
&\big[\,st\big(\, \M^{\tilde{V}_\lam,*}_{\chi,n,(1^d),(1^d)}(U_\lam,dS)\, \cap \, B_\lam\,\big)\,\big] \notag \\
=\
&\sum_s\,\frac{|s|}{\ell(s)!}\,(\xi_{\ell(s)})_*\,
\big[\,st\big(\, (ev_{s,U_0}^*)^{-1}(\triangle_s)\, \cap\, B_0\,\big)\,\big] \ \in\ H_*(\TM_{\chi,n+2d})
\end{align}
where the sum is over all $\ell(s)=d$. Here,
the factor $|s|$ is the degree of the covering map $\pi_{s,\lam}$,
the factor $\ell(s)!$ reflects the fact that each map in the space $\M_s^{\delta,*}(U_\lam)$
in the diagram (\ref{smoothing-local}) is a  labeled map as in (\ref{action}) and
the classes $[st(\,\cdot\,)]$ are defined by the images of cut-down moduli spaces under stabilization map
as in Remark~\ref{cut-down-class}.

In general, if there are maps in the space (\ref{msp-last}) that are not $\delta$-flat then
there is a correction term in (\ref{Thm10.4}) given by three $S$-matrices in $\F_0$ (cf. Definition 11.3 of \cite{IP2})
 for $V_1,V_0$ and $V_2$.
In our case, since the constraint $(F^*)^n$ is supported off the neck,
Lemma 14.6 of \cite{IP2} and Theorem 12.3 of \cite {IP2} imply that the correction term is trivial.
Consequently, by the splitting of the diagonal $\triangle_s$, together with
(\ref{pt-constraint}), (\ref{ILL-ruled-0}) and Remark~\ref{partition},
it follows from (\ref{Thm10.4}) that
\begin{align}\label{sum-local-1}
&GT_{(1^d),(1^d)}^{loc,h,p}\big(\,\prod_{i=1}^n \phi_i^{k_i}(F^*)\,\big) \,=\,
GT_{(1^d),(1^d),\chi}^{loc,h,p}\big(\,(F^*)^n\,\big)\ \cap\  \prod_{i=1}^n \phi_i^{k_i}(F^*)  \notag\\
=\
&\sum\,\frac{|m|}{m!}\,
GT_{(1^d),m,\chi_1}^{loc,h,p}\big(\,(F^*)^{n_1}\,\big)\otimes
GT_{m,(1^d),\chi_2}^{\F_0}\big(\,(F^*)^{n_2}\,\big)\ \cap\ (\xi_{\ell(m)})^*\big(\prod_{i=1}^n \phi_i^{k_i}\big)
\end{align}
where the  sum is over all partitions $m$ of $d$ and $\chi_1+\chi_2-2\ell(m)=\chi$
(cf. Theorem 12.3 of \cite{IP2}).

\bigskip
\non
{\bf Proof of Theorem~A :}
Let $f=(f_1,f_2)$ be a map that contributes to the right hand-side of (\ref{sum-local-1}). Then we have
\begin{itemize}
\item
every connected component of the domain of $(V_0\sqcup V_i)$-regular map $f_i$ has at least two contact (marked) points with $V_0\sqcup V_i$,
\item
every connected component  of the domain  of  $f_2$ has  exactly one contact marked point with $V_0$;
the contact constraint with $V_0$ is $C_{{pt}^{\ell(m)}}$ (see (\ref{rel-F0})) and
the image of each connected component can't pass through more than two distinct points on $V_0$.
\end{itemize}
Noting the gluing map $\xi_{\ell(m)}$ is the map obtained by successively applying gluing maps
as in (\ref{gluingmap}) to connected components,  by Lemma~\ref{L:AC}\,(b) we have
\begin{equation}\label{split-phi}
(\xi_{\ell(m)})^*\big(\prod_{i=1}^n \phi_i^{k_i} \big)\
=\
\prod_{i=1}^{n_1}\phi_i^{k_i}\,\otimes \,
\prod_{i=1}^{n_2} \phi_i^{k_{n_1+i}}.
\end{equation}
The sum formula (\ref{sum-local-1}) together with (\ref{split-phi}) gives
\begin{equation}\label{oursum}
GT_{(1^d),(1^d)}^{loc,h,p}\big(\,\prod_{i=1}^n \phi_i^{k_i}(F^*)\,\big)
\,=\,
\sum_{m} \,\frac{|m|}{m!}\,
GT_{(1^d),m}^{loc,h,p}\big(\,\prod_{i=1}^{n_1} \phi_i^{k_i}(F^*)\,\big)
\cdot
GT_{m,(1^d)}^{\F_0}\big(\,\prod_{i=1}^{n_2} \phi_i^{k_{n_1+i}}(F^*)\,\big).
\end{equation}
Now, Theorem~A  follows from (\ref{oursum}) and Proposition~\ref{dec=rel} in the next section.
\qed

\vskip 1cm

\setcounter{equation}{0}
\section{Descendent Invariants vs. Relative Invariants}
\label{S:D-P}
\bigskip

The aim of this section is to show\,:

\begin{prop}\label{dec=rel}\ \
${\displaystyle
GT_d^{loc,h,p}\big(\,\prod_{i=1}^n\tau_{k_i}(F^*)\,\big)\,=\,
\frac{1}{(d!)^2}\,GT_{(1^d),(1^d)}^{loc,h,p}\big(\,\prod_{i=1}^n\phi_i^{k_i}(F^*)\,\big).
}$
\end{prop}

\bigskip

By the relation of GT and GW invariants (cf. Section 2 of \cite{IP2}),
it suffices to prove Proposition~\ref{dec=rel} for local GW invariants that count maps with connected domains.
Let
$GW_d^{loc,h,p}(\,\cdot\,)$ and  $GW_{m^1,m^2}^{loc,h,p}(\,\cdot\,)$
denote absolute and relative local GW invariants and let
\begin{equation}\label{2dforget}
\pi\,=\,\pi_{2d}\, :\, \CM_{g,n+2d}\,\to\, \CM_{g,n}
\end{equation}
be the forgetful map that forgets the last $2d$ marked points.

\begin{lemma}\label{step1h=0}
If $\sum k_i=0$ or $n\geq 3$, then we have
\begin{equation}\label{E:step1h=0}
GW_{d}^{loc,h,p}\big(\,\prod_{i=1}^n \phi_i^{k_i}(F^*)\,\big)\,=\,
\frac{1}{(d!)^2}\,GW_{(1^d),(1^d)}^{loc,h,p}\big(\,\prod_{i=1}^n \pi^*\phi_i^{k_i}(F^*)\,\big).
\end{equation}
\end{lemma}

\pf
Consider the symplectic fiber sum $\P_h=\F_0\#_{V_1}\P_h\#_{V_2}\F_0$ where
$V_1$ and $V_2$ are two distinct fibers of $\P_h$. This sum can be obtained by
blowing up $\P_h\times \cx$  along $(V_1\cup V_2)\times \{0\}$.
The same arguments of Section~\ref{S:Blow-Up} thus give a sum formula
that is analogous to (\ref{sum-local-1})\,:
\begin{align}\label{lastsum-R}
&GT_{d}^{loc,h,p}\big(\,\prod_{i=1}^n \phi_i^{k_i}(F^*)\,\big)\notag \\
=\
&\sum\,\frac{|m^1||m^2|}{m^1!\,m^2!}\,
GT_{m^1,\chi_1}^{\F_0}\otimes
GT_{m^1,m^2,\chi_0}^{loc,h,p}\big(\,(F^*)^{n}\,\big)\otimes
GT_{m^2,\chi_2}^{\F_0}\ \cap\ (\xi_{\ell(m^1),\ell(m^2)})^*\big(\prod_{i=1}^n \phi_i^{k_i}\big)
\end{align}
where the  sum is over  all $\chi_1+\chi_0+ \chi_2-2\ell(m^1)-2\ell(m^2)=2d(1-h)-2\sum k_i$ and
$\xi_{\ell(m^1),\ell(m^2)}$ is the gluing map obtained by identifying contact points of domains
(see above Theorem~\ref{T:smoothing}).
If $\sum k_i=0$ then for $k=1,2$, by dimension count, we have
\begin{equation}\label{dimcount}
0\,=\,2d-\tfrac12 \chi_k +(\ell(m^k)-d)-\ell(m^k)\,=\,d-\tfrac12 \chi_k.
\end{equation}
This shows $\chi_k=2d$ and hence  $m^k=(1^d)$. Thus, (\ref{E:step1h=0}) for $n=0$ follows from (\ref{com-F0}) and
(\ref{lastsum-R}).

Assume $n\geq 3$ and let $f=(f_1,f_0,f_2)$ be a map that contributes to the right hand side of  (\ref{lastsum-R}).
In order to obtain a sum formula for local GW invariants, we assume that the domain of $f$ is connected.
We have
\begin{itemize}
\item
since all marked points of $f$ map into the middle $\P_h$ side, the domain of $f_k$ ($k=1,2$) mapped into $\F_0$
 has no marked points except contact points,
\item
as in the proof of Theorem~A,
every connected components of the domain of $f_k$
has  one contact point with $V_k$.
\end{itemize}
It follows that the domain of $f_0$ is connected and
the gluing map $\xi=\xi_{\ell(m^1),\ell(m^2)}$ can be obtained by composing gluing maps
as in (\ref{gluingmap}) with $n_2=0$\,:
$$
\eta\,: \CM_{g_1,n_1+1}\times \CM_{g_2,1}\ \to\ \CM_{g_1+g_2,n_1+1}
$$
where $n_1\geq n\geq 3$.
By Lemma~\ref{L:AC}\,(a) for $g_2=0$ and by Lemma~\ref{L:AC}\,(b) for $g_2\geq 1$,
one can see that the pull-back class $\xi^*\phi_i$ restricts to the  trivial class on  two $\F_0$ sides
and hence
$f_k$ is constrained by only $\ell(m^k)$ point contact constraints.
The dimension count (\ref{dimcount}) then shows
$\chi_k=2d$ and $m^k=(1^d)$, and hence $\xi$ is a composition of gluing maps $\eta$ as above with
$(g_2,n_2)=(0,0)$.
Consequently, again by Lemma~\ref{L:AC}\,(a), we have
\begin{equation}\label{pullbackclassbyxi}
\xi^*\phi_i\, =\,1\otimes \big(\,\phi_i\,-\,\sum \delta_{\{i\}\cup I}\,)\otimes 1
\end{equation}
where the sum is over all $I\subset\{n+1,\cdots,n+2d\}$ with $I\ne \emptyset$.
On the other hand, (\ref{AC1}) shows
\begin{equation}\label{pullbackclassbyxi-1}
\pi^*\phi_i\,=\,\phi_i\,-\,\sum \delta_{\{i\}\cup I}.
\end{equation}
Now, (\ref{step1h=0}) for $n\geq 3$
follows from (\ref{com-F0}), (\ref{lastsum-R}), (\ref{pullbackclassbyxi}) and (\ref{pullbackclassbyxi-1}).
\qed

\medskip
\begin{rem}\label{lastremark}
The same argument of the proof of Lemma~\ref{step1h=0}
applies to various sum formulas for the fiber sum of $\P_h$ and $\F_0$.
In particular, for dimension zero local invariants (i.e. $\sum k_i=0$), one can use
(\ref{com-F0}) and the dimension count (\ref{dimcount}) to show
\begin{equation*}
GT_d^{loc,h,p}\,=\,\frac{1}{d!}\,GT_{(1^d)}^{loc,h,p}\,=\,\frac{1}{(d!)^2}GT_{(1^d),(1^d)}^{loc,h,p}
\ \ \ \ \ \mbox{and}\ \ \ \ \
GT_m^{loc,h,p}\,=\,\frac{1}{d!}\,GT_{m,(1^d)}^{loc,h,p}
\end{equation*}
\end{rem}

\medskip
\begin{lemma}\label{step2R}
If $\sum k_i=0$ or $n\geq 3$, then we have
\begin{equation}\label{E;step2R}
GW_{d}^{loc,h,p}\big(\,\prod_{i=1}^n \tau_{k_i}(F^*)\,\big)\,=\,
\frac{1}{(d!)^2}\,GW_{(1^d),(1^d)}^{loc,h,p}\big(\,\prod_{i=1}^n \phi_i^{k_i}(F^*)\,\big).
\end{equation}
\end{lemma}

\pf
If $\sum k_i=0$ then (\ref{E;step2R}) follows from Lemma~\ref{step1h=0}.
Assume  $n\geq 3$ and $h>0$ (we will give  a proof for the case when $h=0$  in the appendix).
Let $V=V_1\cup V_2$ be a union of two distinct fibers of $\P_h$ and denote by
$$
\M^V\,=\,\M^V_{g,n, (1^d),(1^d)}(U,dS)
$$
the local relative GW moduli space.
Let $B$ be  a product of $n$ generic fibers of $\P_h$ each of which is disjoint with $V$ and
let $f$ be a limit map of a sequence in the cut-down moduli space $C\M^{V}\cap B$
where $C\M^V$ is the closure of $\M^V$ in $\CM_{g,n+2d}(U,dS)$.
Then, Remark~\ref{h>0} shows that every  genus zero irreducible component of $f$
maps entirely into either $B$ or $V$. This implies that for $1\leq i\leq n$ and for any
$I\subset\{n+1,\cdots,n+2d\}$ with $I\ne \emptyset$ we have
\begin{equation}\label{sss}
\big[\,\M^V\,\big]\ \cap \  \delta_{\{i\}\cup I} \otimes (F^*)^n\,=\,0.
\end{equation}
Therefore, (\ref{E;step2R}) for $h>0$ follows from (\ref{t=p}), Lemma~\ref{step1h=0}, (\ref{pullbackclassbyxi-1}) and
and (\ref{sss}). \qed

\bigskip
\non
{\bf Proof of Proposition~\ref{dec=rel} :}
We will show that (\ref{E;step2R}) holds for all $n$.
By (\ref{oursum}), we have
\begin{align}\label{lastsum}
GW_{(1^d),(1^d)}^{loc,h,p}\big(\,\prod_{i=1}^n \phi_i^{k_i}(F^*)F^*F^*\,\big)\,&=\,
\frac{1}{d!}\,GW_{(1^d),(1^d)}^{loc,h,p}\big(\,\prod_{i=1}^n \phi_i^{k_i}(F^*)\,\big)\cdot
GT_{(1^d),(1^d)}^{\F_0}\big(\,F^*F^*\,\big) \notag\\
&=\,
d^2\,GW_{(1^d),(1^d)}^{loc,h,p}\big(\,\prod_{i=1}^n \phi_i^{k_i}(F^*)\,\big)
\end{align}
where the second equality follows from Divisor Axiom and (\ref{com-F0}).
On the other hand, one can see that the generalized Divisor Axiom (cf. Lemma 1.4 of \cite{KM}) for descendant GW invariants also
holds for descendant local invariants. Thus, we have
\begin{align*}
&GW_{d}^{loc,h,p}\big(\,\prod_{i=1}^n \tau_{k_i}(F^*)\,\big)\,=\,
\frac{1}{d^2}\,GW_{d}^{loc,h,p}\big(\,\prod_{i=1}^n \tau_{k_i}(F^*)F^*F^*\,\big)\notag \\
=\
&\frac{1}{d^2(d!)^2}\,GW_{(1^d),(1^d)}^{loc,h,p}\big(\,\prod_{i=1}^n \phi_i^{k_i}(F^*)F^*F^*\,\big)\,=\,
\frac{1}{(d!)^2}\,GW_{(1^d),(1^d)}^{loc,h,p}\big(\,\prod_{i=1}^n \phi_i^{k_i}(F^*)\,\big).
\end{align*}
where the second equality follows from Lemma~\ref{step2R} and the last from (\ref{lastsum}).
This completes the proof of Proposition~\ref{ap}. \qed

\vskip 1cm

\setcounter{equation}{0}
\section{Local Contributions to GT invariants of  Ruled Surfaces}
\label{S:Ruled-II}
\bigskip

Let $ \pi:\F_1=\P(\O_E(1)\oplus\O_E)\to E$ be a ruled surface over $E=\P^1$.
This section describes local contributions to GW invariants of $\F_1$ that
is needed for the proof of Theorem~B.
We also denote by $E$ the zero section of $\F_1$ and by $S$ the section class represented by the zero section $E$.

\begin{rem}\label{Renormal}
Let $f:C\to \F_1$ be a holomorphic map
from a smooth domain $C$  that represents the class $dS$.
Then $f$ defines a holomorphic section $\xi$ of the line bundle $(\pi\circ f)^*\O_E(1)$
over $C$.
The zero set $Z(\xi)$ of  $\xi$ is the preimage $f^{-1}(E)$ of the zero section $E$, so
if the image of $f$ does not lie in $E$ (i.e. $\xi\ne 0$)  then
$$
|f^{-1}(E)|\,=\,|Z(\xi)|\,\leq\, \#Z(\xi)\,=\,\deg((\pi\circ f)^*\O_E(1))\,=\,d
$$
where $\#Z(\xi)$ is the number of zeros of $\xi$ counted with multiplicities.
\end{rem}

\medskip

Choose distinct fibers $V_1$, $V_0$ and $V_2$ of $\F_1$ and set
$$
V\,=\,V_1\,\sqcup\, V_0\,\sqcup V_2.
$$

\begin{rem}\label{composition}
Let $f:\P^1\to \F_1$ be a holomorphic map
that represents the class $2S$. If $f$ has a contact vector $(2)$ with $V_i$ at $p_i\in\P^1$ then
for the composition map
\begin{equation*}\label{composition}
\P^1\ \overset{f}{\longrightarrow}\ \F_1\ \overset{\pi}{\longrightarrow}\ E=\P^1
\end{equation*}
the point $p_i$ is a ramification point of multiplicity two.
Thus $f^{-1}(V)$ consists of at least four  points. For otherwise,
$\pi\circ f$ is a holomorphic map of degree 2 with three ramification points,
which is impossible by  the Riemann-Hurwitz formula.
\end{rem}

Remark~\ref{Renormal} and Remark~\ref{composition} immediately give\,:

\begin{lemma}\label{IL-pc}
Let $f$ be a holomorphic map into $\F_1$
from a smooth domain of genus $g$ representing the class $dS$ and satisfying
$$
\im(f)\cap (V_1\sqcup V_2) \,=\, (V_1\sqcup V_2)\cap E.
$$
If either (i) $d=1$ or (ii) $d=2$, $g=0$ and the contact vector of $f$ with $V_0$ is $(2)$
then the image of $f$  lies in $E$.
\end{lemma}

\medskip

Fix a neighborhood $U$ of  the zero section $E$  whose closure is disjoint from the infinity section of $\F_1$.
For $V$-compatible $(J,\nu)$ and ordered sequences $s^1$ and $s^2$ with $\deg(s^1)=\deg(s^2)=2$,
choose $\ell(s^1)$ points $\{p_i\}$ in $V_1$ and $\ell(s^2)$ points $\{q_j\}$ in $V_2$  and
let
\begin{equation}\label{cut-down}
\M_{0,s^1,(2),s^2}^{V_1,V_0,V_2}(U,2S,J,\nu)\,\cap\,\{p_i,q_j\}\ \subset\
\M_{0,s^1,(2),s^2}^{V_1,V_0,V_2}(U,2S,J,\nu)
\end{equation}
denote the cut-down moduli space of  $V$-regular maps $f$ into $U$
with (connected) domain of genus zero and with $n$ contact points $\{x_i\}$ satisfying\,:
$f(x_0)\in V_0$, $f(x_i)=p_i$ and $f(x_{\ell(s^1)+j})=q_j$
where $n=\ell(s^1)+\ell(s^2)+1$.
This cut-down moduli space  has
(formal) dimension zero. Consider a sequence of  maps $(f_k,C_k;\{x_i^k\})$ in the cut-down moduli spaces
\begin{equation}\label{cut-down-seq}
\M_{0,s^1,(2),s^2}^{V_1,V_0,V_2}(U,2S,J_k,\nu_k)\,\cap\,\{p_i^k,q_j^k\}
\end{equation}
where the points $\{p_i^k\}\sqcup \{q_j^k\} \subset (V_1\sqcup V_2)$ converge
to  points in $(V_1\sqcup V_2)\cap E$ and
$V$-compatible  $(J_k,\nu_k)$ converges to the complex
structure of $\F_1$  as $k\to \infty$.
By the Gromov Compactness Theorem, after passing to subsequences, the sequence of maps
$(f_k,C_k;\{x_i^k\})$ then converges to
a holomorphic map
\begin{equation}\label{limit-map}
(f,C;\{x_i\})\in\CM_{0,n}(\overline{U},2S).
\end{equation}

Since $\im(f)\subset \overline{U}$, every component of $C$ mapped entirely into
$V$ is a  ghost component.
The following lemma shows that the image of $f$ lies in the zero section $E$.
Let $s^0=(2)$.

\begin{lemma}\label{IL-cut-down}
Let $(f,C;\{x_i\})$ be  as above and let  $C_i$ denote an irreducible component $C$
that  contains a marked point $x_{i_0}$ mapped into $V_i$.
If $s^i=(2)$ for some $0\leq i\leq 2$ then we have
\begin{itemize}
\item[(a)]
if $f$ is $V_i$-regular then the restriction $f$ to $C_i$  represents the class $2S$,
\item[(b)]
if $f$ is not $V_i$-regular then
$C_i$ is a ghost component with $x_j\notin C_i$ for $j\ne i_0$ and
$C\setminus C_i$ has two connected components $C_i^\ell$ such that the restriction of $f$ to $C_i^\ell$
represents the class $S$.
\end{itemize}
In particular, the image of $f$ lies in the zero section $E$ of $\F_1$.
\end{lemma}

\pf
First note that, since $\im(f_k)\to \im(f)$ and $\im(f)\subset \overline{U}$, we have
\begin{itemize}
\item[(i)]
$\im(f)\cap(V_1\sqcup V_2) = (V_1\sqcup V_2) \cap E$,
\item[(ii)]
the restriction of $f$ to each component of $C$ represents the class $kS$ where $0\leq k\leq 2$.
\end{itemize}
If a marked point $x_i$ is a limit point of the contact points $x_i^k$ of $f_k$ with $V_i$, i.e.
$f_k(x_i^k)\in V_i$, then $x_i$ is a contact point of $f$ with the same contact order as $x_i^k$
unless a component containing $x_i$ maps into $V$.
Thus, if $f$ is $V_i$-regular then $f$ has a contact vector (2)
with $V_i$ at $x_{i_0}$ and hence by (ii) the restriction of $f$ to $C_i$
represents the class $2S$. This proves (a).

On the other hand, if $f$ is not $V_i$-regular then
$C_i$ is a ghost component mapped into $V_i$. Moreover,
by the assumption $s^i=(2)$ only one marked point $x_{i_0}$ maps into $V_i$, so
$x_j\notin C_i$ for $j\ne i_0$.
Note that, since $C$ is a connected curve of (arithmetic) genus zero,
each irreducible component of $C$ is smooth of genus zero.
Thus Remark~\ref{ghost} implies that $C\setminus C_i$ has  at least two connected components.
Let $C_i^\ell$ be a connected component of $C\setminus C_i$. Then,
$$
|C_i^\ell \cap \overline{C\setminus C_i^\ell}|\,=\,|C_i^\ell \cap C_i|\,=\,1
$$
where the second equality follows from the fact $C$ is a connected curve of genus zero.
So, if $C_i^\ell$ maps to a point then $f(C_i^\ell)\in V_i$ and, by Remark~\ref{ghost},
$C_i^\ell$ has at least two marked points $x_j$ with $f(x_j)\in V_i$.
This is impossible since $f(x_j)\in V_j$ for some $V_j$ disjoint with $V_i$.
Thus there are two connected components of $C\setminus C_i$ such that the restriction of $f$ to both components
represent the class $S$. This completes the proof (b).

Now (a), (b), (i), (ii) and Lemma~\ref{IL-pc} imply that the image of $f$ lies in the zero section
$E$ since $s^0=(2)$. \qed

\begin{rem}\label{rem-cal}
Let $f$ be a limit map as in Lemma~\ref{IL-cut-down} and suppose that
$C_0$ and $C_1$ are  irreducible components of the domain $C$ of $f$
which contain  marked points $x_0$ and $x_1$ mapped into $V_0$ and $V_1$.
Suppose $s^1=(2)$. Then, since $s^0=(2)$, Lemma~\ref{IL-cut-down} implies that
$f$ is $V_0$-regular if and only if $f$ is $V_1$-regular.
Suppose $f$ is not $(V_0\sqcup V_1)$-regular. Then, by Lemma~\ref{IL-cut-down}\,(b),
both $C_0$ and $C_1$ are ghost components with $C_0\cap C_1=\emptyset$.
Since the domain $C$ is a connected curve of genus zero,
the same argument for the proof of Lemma~\ref{IL-cut-down}\,(b) shows that
$C\setminus (C_0\sqcup C_1)$ has at least three connected components $C^\ell_{01}$ such that
the restriction of $f$ to each $C^\ell_{01}$ represents the class $S$.
This is impossible, so (i) $f$ must be  $(V_0\sqcup V_1)$-regular  (ii) $C_0=C_1$ and
(iii)  the restriction of
$f$ to the component $C_0=C_1$ represents the class $2S$.
In particular, if $s^1=s^2=(2)$ then $f^{-1}(V)$
consists of three contact points. This contradicts  Remark~\ref{composition}.
Therefore, if $s^1=s^2=(2)$  then for all large $k$ the cut-down moduli space (\ref{cut-down-seq}) must be empty.
\end{rem}

\medskip
For $V$-compatible $(J,\nu)$ and for points $\{p_i\}$ in $V_1$ and $\{q_j\}$ in $V_2$ denote by
$$
\M \,=\,\M_{0,s^1,(2),s^2}^{V_1,V_0,V_2}(\F_1,2S)\,\cap\,\{p_i,q_j\}
$$
the cut-down moduli space of $V$-regular $(J,\nu)$-holomorphic maps  into $\F_1$ with the point constraints
$\{p_i,q_j\}$. Then we have a splitting of contributions to relative GW invariant
\begin{equation}\label{decomp-contrib}
GW^{V_1,V_0,V_2}_{\F_1,2S, 0,s_1,(2),s_2}(C_{pt^{\ell(s^1)}},C_F,C_{pt^{\ell(s^2)}})\,=\,
[\,\M\,]\,=\,[\,\M(U)\,]\,+\,[\,\M\setminus \M(U)\,]
\end{equation}
where $\M(U)$ is the cut-down moduli space (\ref{cut-down}) and
$\M\setminus \M(U)$ consists of maps $f$ in $\M$ whose image dose not lie in $U\subset \F_1$.
Here, $[\,\cdot\,]$ denotes the zero-dimensional homology class defined by cut-down moduli space as in
Remark~\ref{cut-down-class}.
In general, this splitting is not well-defined, namely it depends on the choice of $(J,\nu)$ and the point constraints
$\{p_i,q_j\}$. However, by the Gromov Compactness theorem and Lemma~\ref{IL-cut-down},
the splitting (\ref{decomp-contrib}) is well-defined whenever
(i) $(J,\nu)$ is   close to the complex structure
of $\F_1$ and (ii) points $\{p_i,q_j\}$ are also close to
points in $(V_1\sqcup V_2)\cap E$. In such cases, the contribution
\begin{equation}\label{local-F1}
[\,\M(U)\,]\,=\,\big[\,\M_{0,s^1,(2),s^2}^{V_1,V_0,V_2}\big(U,2S;C_{pt^{\ell(s^1)}},C_F,C_{pt^{\ell(s^2)}}\big)\,\big]
\end{equation}
is  independent of the choice of $(J,\nu)$ and $\{p_i,q_j\}$.
Furthermore, this local contribution is independent of the choice of the neighborhood $U$ of the zero section $E$ in $\F_1$.

\begin{lemma}\label{comp-ruled}
For any neighborhood $U$ of the zero section $E$ in $\F_1$ whose closure is disjoint from the infinity section of $\F_1$,
we have
$$
(a)\ \ \big[\,\M_{0,(2),(2),(2)}^{V_1,V_0,V_2}\big(U,2S;C_{pt},C_F,C_{pt}\big)\,\big]=0\ \
\ \ \ \
(b)\ \ \big[\,\M_{0,(2),(2),(1,1)}^{V_1,V_0,V_2}\big(U,2S;C_{pt},C_F,C_{pt^2}\big)\,\big]=1
$$

\end{lemma}

\pf (a) follows from Remark~\ref{rem-cal}. Let $s^1=(2)$, $s^2=(1,1)$ and $(f,C;\{x_i\})$ be a limit map as in (\ref{limit-map}).
Again by Remark~\ref{rem-cal}, $f$ is $(V_1\sqcup V_0)$-regular map into $E\subset \F_1$ and
the restriction of $f$ to the component $C_0=C_1$ containing  $x_0$ and $x_1$
represents the class $2S$. Stability of $f$ then implies either  $C=C_0$  or
$C=C_0\cup C_2$ where
$C_2$ is a ghost component containing the marked points $x_2$ and $x_3$.
Suppose $C=C_0\cup C_2$. In this case, $C$ has one node since the (arithmetic) genus of $C$ is zero.
The restriction $f_0=f_{|C_0}$ has  a contact order two with $V_2$ at the node of $C$,
so $f_0^{-1}(V)$ consists of three contact points. This contradicts Remark~\ref{composition}.
Therefore,
$f$ is a holomorphic map from $C=\P^1$ into $E=\P^1$  of degree two with two ramification points $x_0$ and $x_1$
and $f(x_i)\in V_2\cap E$ for $i=2,3$. Observe that there is a unique such map $f$.

Let  $J$ denote the complex structure on $\F_1$ and set
$$
\M^V \!= \M^{V_1,V_0,V_2}_{0,(2),(2),(1,1)}(U,2S,J)
\ \ \ \ \mbox{and}\ \ \ \
h=h_{(2)}\times h_{(1,1)}:\M^V\ \to\  V_1\times (V_2\times V_2)
$$
where $h_{(2)}$ and $h_{(1,1)}$ are
evaluation maps as in (\ref{ev-h}).
Let $D_f$ be the (full) linearization  of holomorphic map equation at $f$.
Since the normal bundle of $\im(f)=E$ is $\O_E(1)$,
we have
\begin{itemize}
\item
$\cok D_f=H^1(f^*\O_E(1))=H^1(\O_{\P^1}(2))=0$ and hence
\item
$\M^V$ is smooth near $f$ with $T_f\M^V=H^0(f^*\O_E(1))=H^0(\O_{\P^1}(2))$,
\item
$dh_f(\xi)=(\,\xi(x_1),\xi(x_2),\xi(x_3)\,)$.
\end{itemize}
In fact, regarding the neighborhood $U$ of $E$ in $\F_1$ as a
disk subbundle of $\O_E(1)$, one can identify holomorphic sections
$\xi$ of $f^*U\subset f^*\O_E(1)$ with $V$-regular holomorphic maps $f_\xi$ in  $\M^V$ ---
in local trivialization,   $f_\xi(x)=(f(x),\xi(x))$.
By Remark~\ref{Renormal} the differential $dh_f$ is one-to-one. Thus $dh_f$  is  onto since $h^0(\O_{\P^1}(2))=3$.
We can now conclude that the contribution of $(f,\P^1;\{x_i\})$ to the invariant (\ref{local-F1}) is $+1$ since
$f$ has no nontrivial automorphisms  and
$dh_f$ is onto and complex linear.
This completes the proof of (b).
\qed

\medskip
We will compute the local contribution (\ref{local-F1}) for the case $s^1=s^2=(1,1)$ in the next section
(see (\ref{com-by-rel})).

\setcounter{equation}{0}
\section{Spin Curve Degeneration and Sum Formula I}
\label{S:SC-Sum}
\bigskip

This section proves  Theorem~B\,(a) in three steps. First we review the {\em dualizing sheaf}.
The dualizing sheaf $\omega_X$ of a variety $X$ (if exists)
is the unique invertible sheaf making Serre duality valid;
when $X$ is smooth $\omega_X$ is the canonical bundle $K_X$.
For a proper holomorphic map $f:X\to B$ between two smooth varieties,
the relative dualizing sheaf $\omega_f$ is the locally free
rank one sheaf $\omega_X\otimes (f^*\omega_B)^{-1}$ whose restriction to each fiber $X_b$ is the dualizing
sheaf of $X_b$ (cf. \cite{T} and \cite{HM}).

\bigskip

\non
{\bf Step 1 :}
Let $D_1\cup_q D_2$ be a union of two smooth curves $D_1$ and $D_2$ of genera $h_1$ and $h_2$,
meeting at one point $q$. Blowing up the point $q$  yields a nodal curve $D_0$ with an exceptional component
$E=\P^1$ that meets $\bar{D}=D_1\sqcup D_2$ at two points.
A theta characteristic of the nodal curve $D_0$ is a line bundle
$N_0$ together with a homomorphism $\phi:N_0^2\to \omega_{D_0}$ satisfying\,:
\begin{itemize}
\item
$N_0$ restricts to $\O(1)$ on the exceptional component $E$,
\item
$\phi$ vanishes identically on  $E$ and restricts to an  isomorphism
$N_0^2|_{\bar{D}} \simeq \omega_{\bar{D}}$,
\end{itemize}
where $\omega_{D_0}$ and $\omega_{\bar{D}}$ are the dualizing sheaves of $D_0$ and $\bar{D}$ respectively.
Since $\omega_{\bar{D}}$ is the canonical bundle of $\bar{D}$, the line bundle $N_0$ restricts to
theta characteristic $N_1$ on $D_1$ and restricts to theta characteristic $N_2$ on $D_2$.
The triple $(D_0, N_0,\phi)$ is  a {\em spin curve} of genus $h=h_1+h_2$ with parity $p\equiv p_1 + p_2$ (mod 2) where
$p_i$ is a parity of the spin curve $(D_i,N_i)$.
It then follows from a universal deformation of the spin curve $(D_0,N_0,\phi)$ (cf. pg 570 \cite{C}) that there are
\begin{itemize}
\item
a family of curves $\rho:{\cal D}\to \Delta$ where $\Delta$ is a unit disk in $\cx$,
the fiber $D_\lam$ over $\lam\ne 0$ is a smooth curve of genus $h$ and the central fiber is the nodal curve
$D_0$,
\item
a line bundle $\pi:{\cal N}\to {\cal D}$ together with a homomorphism ${\Phi}: {\cal N}^2\to \omega_\rho$
such that each
$\big({\cal D}_{|\rho^{-1}(\lam)},{\cal N}_{|\rho^{-1}(\lam)}, \Phi_{|\rho^{-1}(\lam)}\big)$
is a spin curve of genus $h$ with parity $p$,
\end{itemize}
where $\omega_\rho$ is the relative dualizing sheaf of $\rho$ (for more details see \cite{C}).

Let  ${\cal N}_{D_\lam}$ ($\lam\ne 0$) be the total space of ${\cal N}_{\,|D_\lam}$ and
${\cal N}_{\bar{D}}$ be the total space of ${\cal N}_{\,|\bar{D}}$.
Since both total spaces ${\cal N}_{D_\lam}$ and ${\cal N}_{\bar{D}}$ are smooth, there are short exact sequences
$$
0\to \pi^*{\cal N}_{\,|D_\lam}\to T{\cal N}_{D_\lam}\to \pi^*TD_\lam\to 0
\ \ \ \ \ \mbox{and}\ \ \ \ \
0\to \pi^*{\cal N}_{\,|\bar{D}}\to T{\cal N}_{\bar{D}}\to \pi^*T\bar{D}\to 0.
$$
It then follows from these exact sequences that
\begin{equation}\label{can-formulas}
K_{{\cal N}_{D_\lam}} = \pi^*{\cal N}_{\,|D_\lam}^*\otimes \pi^*K_{D_\lam}
\ \ \ \ \  \mbox{and}\ \ \ \ \
K_{{\cal N}_{\bar{D}}} = \pi^*{\cal N}_{\,|\bar{D}}^*\otimes \pi^*K_{\bar{D}}.
\end{equation}
On the other hand, the homomorphism
$\Phi:{\cal N}^2\to \omega_\rho$ induces a homomorphism
$$
\Phi^\prime : \pi^*{\cal N} = \pi^*({\cal N}^*\otimes {\cal N}^2) \ \longrightarrow\
\pi^*({\cal N}^*\otimes \omega_\rho).
$$
Let $\sigma$ be the tautological section of $\pi^*{\cal N}$.
Since the relative dualizing sheaf $\omega_\rho$ restricts to dualizing sheaf $\omega_{D_\lam}$ on  $D_\lam$,
by (\ref{can-formulas})  the composition  $\Phi^\prime\circ \sigma$ is a section of
$\pi^*({\cal N}^*\otimes \omega_\rho)$ satisfying\,:
\begin{itemize}
\item
$\Phi^\prime\circ \sigma$ restricts to a holomorphic 2-form $\a_\lam$ on
${\cal N}_{D_\lam}$ ($\lam\ne 0$) whose zero set is $D_\lam$,
\item
$\Phi^\prime\circ \sigma$ restricts to a holomorphic 2-form
 $\a_0$ on  ${\cal N}_{\bar{D}}$ whose zero set is $\bar{D}$.
\end{itemize}
Let ${\cal N}_E$ denote the total space of ${\cal N}_{\,|E}$.
For sufficiently small $\epsilon>0$,
choose a bump function $\beta$
that is 1 on the complement of the $2\epsilon$-neighborhood of  ${\cal N}_E$
in ${\cal N}_{\cal D}$ and
vanishes on the $\epsilon$-neighborhood ${\cal U}_\epsilon$ of ${\cal N}_E$.
For each point $x\in {\cal N}_{\cal D}\setminus {\cal U}_\epsilon$ we set
$\beta\a_\lam(v_x,\cdot\,)=0$ for normal vectors $v_x$ to submanifold ${\cal N}_{D_\lam}$
(or ${\cal N}_{\bar{D}}$ for $\lam=0$) at $x$ and
define $\a=\beta\a_\lam$ at $x$. Extension $\a$ by zero then
gives a 2-form on ${\cal N}_{\cal D}$ that restricts to $\b\a_\lam$ on
${\cal N}_{D_\lam}$ ($\lam\ne 0$) and to $\b\a_0$ on ${\cal N}_{\bar{D}}$.

\bigskip
\non
{\bf Step 2 :}
Consider the projectivization $\P({\cal N}\oplus \O_{\cal D})$ over ${\cal D}$ that gives a degeneration
$$
\lam\,:\,Z=\P({\cal N}\oplus \O_{\cal D})\,\to \,{\cal D}\,\to\,\Delta
$$
whose fiber $Z_\lam$ ($\lam\ne 0$) is a ruled surface over $D_\lam$ isomorphic to $\P_h=\P(N\oplus \O_D)$
where $(D,N)$ is a (smooth) spin curve of genus $h=h_1+h_2$ with parity $p\equiv p_1 + p_2$ (mod 2) and
whose central fiber $Z_0$ is the singular (ruled) surface
$$\P_{h_1} \cup_{V_1} \F_1 \cup_{V_2} \P_{h_2}\,\to\,D_0$$
where $V_1$ and $V_2$ are fibers over the nodes of $D_0$.
Note that the general fiber $Z_\lam$ ($\lam\ne 0$) is the symplectic fiber sum
$$
\P_h\,=\,\P_{h_1}\#_{V_1}\#\F_1\#_{V_2}\P_{h_2}.
$$

Let $U$ be an (open) neighborhood of the zero section of $Z=\P({\cal N}\oplus \O_{\cal D})$
and fix an isomorphism $\Psi$ from $U$ to some neighborhood of ${\cal D}\subset {\cal N}_{\cal D}$
taking the zero section of $Z$ to ${\cal D}$.
Choose a point $q_0$ in the exceptional component $E\subset D_0$ that is not a nodal point and
let $B\subset {\cal D}$ be a normal disk to $E\subset {\cal D}$ at $q_0$, namely
the intersection $B\cap D_\lam$ is  one point for all small $|\lam|$.
Let $\tilde{V}_\lam$ be the fiber of $Z_\lam\to D_\lam$ over the intersection point $B\cap D_\lam$
and set
$$ V=\,V_1\,\sqcup\,V_0\,\sqcup\,V_2
$$
where $V_0=\tilde{V}_0$. Denote by $\J(Z)$
the space of all $(J,\nu)$ on $Z$ satisfying\,: (i) each $Z_\lam$ is $J$-invariant and (ii)
the restriction of $(J,\nu)$  to $Z_0$ and $Z_\lam$ ($\lam\ne 0$) are  $V$-compatible
and $\tilde{V}_\lam$-compatible respectively. We will use the same notation $(J,\nu)$ for its restriction to each
$Z_\lam$.
Denote by the same $S$  the section classes of $\P_{h_1}$, $\F_1$, $\P_{h_2}$ and $Z_\lam$ represented by the zero sections.
Let $D_i$ ($i=1,2$) and $D_\lam$ denote the zero sections of $\P_{h_i}$ and
$Z_\lam$ respectively. For each small $|\lam|$ we set
$$
U_\lam\,=\,U\,\cap\, Z_\lam.
$$

By using the  2-form $\a$ on ${\cal N}_D$ together with
the isomorphism $\Psi$, we obtain\,:

\begin{lemma}\label{IL-scd}
There is an almost complex structure $J_{\sv}$ on $Z$ satisfying\,:
\begin{itemize}
\item[(a)]
$(J_{\sv},0)\in\J(Z)$ and $J_{\sv}$ restricts to the complex structure of $\F_1$,
\item[(b)]
$\CM_{\chi,n}^*(U_\lam,dS,J_{\sv})\,=\,\CM_{\chi,n}^*(D_\lam,d)$ \ and\ \ \
$\CM_{\chi,n}^*(\overline{U}\!_0\cap \P_{h_i}, dS, J_{\sv})\,=\,\CM_{\chi,n}^*(D_i,d)$,
\item[(c)]
for generic $(J,\nu)\in \J(Z)$ sufficiently close to $(J_{\sv},0)$ and for small $|\lam|>0$
\begin{equation*}\label{vfc-scd}
\big[\,\M_{\chi,(2)}^{\tilde{V}_\lam,*}(U_\lam,2S)\,\big]\,=\,GT_{(2)}^{loc,h,p}
\ \ \ \ \ \mbox{and}\ \ \ \ \
\big[\,\M_{\chi_i,s_i}^{V_i,*}(\P_{h_i}\cap U_0,2S)\,\big]\,=\,GT_{s_i}^{loc,h_i,p_i }
\end{equation*}
\end{itemize}
where $\chi=2-4h$, $\chi_i=2\ell(s_i)-4h_i$ for $i=1,2$.
\end{lemma}

\pf
By the isomorphism $\Psi$ as above, one can regard $\a_\lam$ ($\lam\ne 0$) and
$\a_0$ as  holomorphic 2-forms on $U_\lam$ and on $U_0\cap (\P_{h_1}\sqcup \P_{h_2})$
whose zero sets are $D_\lam$ and $D_1\sqcup D_2$, respectively. Similarly, one can also regard
$\b$ as a bump function on $U$ and $\a$ as a 2-form on $U$
satisfying (i) $\a$ vanishes on some neighborhoods of $\F_1$, $V$ and $\tilde{V}_\lam$ for small $|\lam|$, and (ii)
the restriction of $\a$ to $U_\lam$ and $U_0\cap (\P_{h_1}\sqcup \P_{h_2})$ are respectively
$\b\a_\lam$ and $\b\a_0$.
Now, let $J_{\sv}$ be the almost complex structure on $Z$ induced by $\a$ and the formula (\ref{Ja}).
Then, (a) follows from (i), (b) follows from (ii) and Remark~\ref{bump-ILL} and (c) follows from definition
and compactness by (b).    \qed

\medskip

\begin{rem}
Let $f$ be a map in $\M_{\chi,(2)}^{\tilde{V}_\lam}(U_\lam,2S)$ with $f^{-1}(\tilde{V}_\lam)=\{x\}$
where the Euler characteristic $\chi=2-4h$.
If there is a  connected component of the domain of $f$ that does not contain the contact point $x$ then
the restriction of $f$ to that component represents the trivial homology class.
The stability of $f$ thus shows  $C$ is connected since
there is no marked points except the contact point $x$. Consequently, for $\lam\ne 0$ we have
\begin{equation}\label{connected-1}
\M_{\chi,(2)}^{\tilde{V}_\lam,*}(U_\lam,2S)\,=\,\M_{\chi,(2)}^{\tilde{V}_\lam}(U_\lam,2S)\,=\,
\M_{g,(2)}^{\tilde{V}_\lam}(U_\lam,2S)
\end{equation}
where the genus $g=2h$. Similar arguments also show that
\begin{equation}\label{connected-2}
\M_{\chi_0,s_1,(2),s_2}^{V_1,V_0,V_2,*}(\F_1,2S)\,=\,
\M_{g_0,s_1,(2),s_2}^{V_1,V_0,V_2}(\F_1,2S)
\end{equation}
where  the Euler characteristic $\chi_0=2$ and the genus $g_0=0$.
\end{rem}

\bigskip
\non
{\bf Step 3 :}
Choose an (open) neighborhood $W$ of the zero section of  $Z=\P({\cal N}\oplus \O_{\cal D})$
with $\overline{W}\subset U$ and for each small $|\lam|$ set
$$   W_\lam=W\cap Z_\lam.   $$
The following is the key fact to the proof of Theorem~B (a).

\begin{lemma}\label{key-scd}
For $(J,\nu)\in \J(Z)$ sufficiently close to $(J_{\sv},0)$ and for small $|\lam|>0$, we have
$$
\M_{\chi,(2)}^{\tilde{V}_\lam}(U_\lam,2S)\,\setminus\,
\M_{\chi,(2)}^{\tilde{V}_\lam}(W_\lam,2S)\,=\,\emptyset
$$
where the Euler characteristic  $\chi=2-4h$.
\end{lemma}

\pf Suppose not. Then, there exists a sequence of $(J_k,\nu_k)$-holomorphic maps $f_k$ into $U_{\lam_k}$
with $\im(f_k)\cap (U_{\lam_k}\setminus W_{\lam_k})\ne \emptyset$
where $\lam_k\to 0$ and $(J_k,\nu_k)$ converges to $(J_{\sv},0)$
as $k\to \infty$.
After passing to subsequences, by the Gromov Compactness Theorem, $f_k$ converges to a $J_{\sv}$-holomorphic
map $f$ into $Z_0$ such that
(i) $\im(f)\subset \overline{U}\!_0$ and (ii) $\im(f)\cap (\overline{U}\!_0\setminus W_0)\ne \emptyset$.
By (i) and Lemma~\ref{IL-scd},
$f$ can be split as $f=(f_1,f_0, f_2)$ where
$f_1$ and $f_2$ are respectively holomorphic maps into $D_1$ and $D_2$, and
$f_0$ is a  holomorphic map into $\F_1$ such that
\begin{equation}\label{intersection}
\im(f_0)\,\cap\, (V_1\sqcup V_2)\,=\,\big(\,\im(f_1)\sqcup \im(f_2)\,\big)\, \cap\, (V_1 \sqcup V_2)
\,=\,E\,\cap\, (V_1 \sqcup V_2).
\end{equation}
Note that the domain $C$ of $f$ is a connected curve of genus $2h$ since
by (\ref{connected-1})
$f$ is a limit of maps with connected domains of genus $2h$.
Also note that if $f$ is not $(V_1\sqcup V_2)$-regular, there is a ghost component mapped into $V_1\sqcup V_2$.

Let $f_{12}=(f_1,f_2)$, $C_{12}$ be the domain of $f_{12}$ and $C_0$ be the domain of $f_0$.
We can assume that $C_0$ contains all ghost components mapped into $V_1\sqcup V_2$.
Then $f_{12}$ is $(V_1\sqcup V_2)$-regular, so $f_{12}^{-1}(V_1\sqcup V_2)=(C_{12}\cap C_0)$.
Let $\ell=|f_{12}^{-1}(V_1\sqcup V_2)|$. Since $g-1=-\frac12 \chi$, we have
\begin{equation}\label{euler-genus}
2h\,=\,g(C)\,=\,-\tfrac12\chi(C_0) -  \tfrac12 \chi(C_{12}) + \ell +1.
\end{equation}
Consider $f_{12}$ as a holomorphic map into $D_1\sqcup D_2$ and
apply the Riemann-Hurwitz formula to each irreducible component of $C_{12}$. This gives
\begin{equation}\label{RH}
-\tfrac12\chi(C_{12}) + \ell \,\geq\, 2h
\end{equation}
since the geometric genus of each irreducible component is less than or equal to its
arithmetic genus and  $f_{12}$ has at least $2(4-\ell)$ ramification points.
Consequently, by (\ref{euler-genus}) and (\ref{RH}) we have
\begin{equation}\label{1-claim}
\chi(C_0)\,\geq \,2.
\end{equation}

Note that the image of $f_0$ does not lie in the zero section $E$ of $\F_1$ by (ii) since $f_{12}$ maps into
$D_1\cup D_2\subset W_0$.
Remark~~\ref{Renormal} and (\ref{intersection})  imply that there is exactly one irreducible component $C_0^\prime$ of
$C_0$ such that the restriction $f_0^\prime=f_{|C_0^\prime}$  has contact vectors (2) with $V_1$ and $V_2$ and
all other irreducible components of $C_0$ are ghost components.
Let $C^1_0$ be the connected component of $C_0$ that contains $C_0^\prime$.
Since $f$ has no degree zero components, we have
\begin{itemize}
\item
all ghost components mapped into $V_0$ are contained in $C_0^1$,
\item
if there exists a connected component $C_0^2\ne C_0^1$ of $C_0$ then  $C_0^2$ is a union of ghost components such that
$C_0^2$ has no marked points and maps into either $V_1$ or $V_2$.
Since $|C_0\cap C_{12}|=\ell\leq 4$ and $C_0^1\cap C_{12}$ contains at least one point mapped into $V_1$ and at least one point
mapped into $V_2$,
we have $|C^2_0\cap \overline{C\setminus C^2_0}|=|C^2_0\cap C_{12}|=1$ and hence $g(C_0^2)>0$ by Remark~\ref{ghost}.
\end{itemize}
Now from (\ref{1-claim}) we have  $g(C_0^1)=0$, so each irreducible component of
$C_0^1$ has genus zero and no two irreducible components meet at more than one point.
In particular, since $C_0^\prime$ has genus zero, Lemma~\ref{IL-pc} implies that
$f_0^\prime$ has a contact vector $(1,1)$ with $V_0$.
In this case, since $f$ is a limit map of a sequence of maps with contact order (2) with $V_0$,
there is a ghost component mapped into $V_0$. Let $C_0^{\prime\prime}$ be a connected component of
the union of all ghost components mapped into $V_0$. Then,
$C^{\prime\prime}_0\subset C_0^1$, so
$g(C^{\prime\prime}_0)=0$ and
$$
|C_0^{\prime\prime}\cap \overline{C\setminus C_0^{\prime\prime} }|\,=\,
|C_0^{\prime\prime}\cap C_0^\prime|\,<\,2.
$$
Since $C_0^{\prime\prime}$ has at most one marked point, we have a contradiction by Remark~\ref{ghost}.
\qed

\bigskip
\medskip
\non
{\bf Proof of Theorem~B\,(a) :}
The proof is identical to the proof of Theorem~A. We only outline the proof.
For  (ordered) sequences $s^i$ with $\deg(s_i)=2$ where $i=1,2$, consider the evaluation map
that records the intersection points with $V_1$ and $V_2$\,:
\begin{align*}
ev_{s^1,s^2,U_0}^* \,: \ &{\textstyle \bigcup}\
\M_{\chi_1,s^1}^{V_1,*}\big(\,\P_{h_1}\cap U_0,2S\,\big)\times
\M_{\chi_0,s^1,(2),s^2}^{V_1,V_0,V_2,*}\big(\,\F_1\cap U_0,2S\,\big)\times
\M_{\chi_2,s^2}^{V_2,*}\big(\,\P_{h_2}\cap U_0,2S\,\big)\\
&\ \longrightarrow\
\Big(V_1^{\ell(s^1)}\times V_1^{\ell(s^1)}\Big )\times
\Big(V_2^{\ell(s^2)}\times V_2^{\ell(s^2)}\Big )^{\phantom{\displaystyle \cup}}
\end{align*}
where the union is over all
$ \chi_1+\chi_0+\chi_2 -2\ell(s^1)-2\ell(s^2)=2-4h.$
Let $\triangle_{s^i}$ be the diagonal
of $V^{\ell(s^i)}\times V^{\ell(s^i)}$ for $i=1,2$.
Lemma~\ref{key-scd} and Theorem~10.1 of \cite{IP2} then give
\begin{equation}\label{sum-scd-2}
\big[\,\M_{\chi,(2)}^{\tilde{V}_\lam,*}(U_\lam,2S)\,\big]\,=\,
\sum_{s^1,s^2}
\frac{|s^1||s^2|}{\ell(s^1)!\,\ell(s^2)!}\,
\big[\,\big(ev_{s^1,s^2,U_0}^*\big)^{-1}(\triangle_{s^1}\times \triangle_{s^2})\,\big]
\end{equation}
where $\chi=2-4h$. 
On the other hand, the splitting of the diagonal $\triangle_{s^1}\times \triangle_{s^2}$ yields
\begin{align}\label{sum-scd-1}
&\sum_{s^1,s^2}
\frac{|s^1||s^2|}{\ell(s^1)!\,\ell(s^2)!}\,
\big[\,\big(ev_{s^1,s^2,U_0}^*\big)^{-1}(\triangle_{s^1}\times \triangle_{s^2})\,\big]\notag \\
=\ &
\sum_{m^1,m^2}\frac{|m^1||m^2|}{m^1!\,m^2!}\,\,
GT_{m^1}^{loc,h_1,p_1}\cdot
\big[\,\M_{0,m^1,(2),m^2}^{V_1,V_0,V_2}(\F_1\cap U_0,2S;C_{pt^{\ell(m^1)}},C_F,C_{pt^{\ell(m^2)}})\,\big]\cdot
GT_{m^2}^{loc,h_2,p_2} \notag \\
=\ &
GT_{(1^2)}^{loc,h_1,p_1}\cdot GT_{(2)}^{loc,h_2,p_2} \,+\,
GT_{(2)}^{loc,h_1,p_1}\cdot GT_{(1^2)}^{loc,h_2,p_2} \notag \\
+&\ \frac14\,
GT_{(1^2)}^{loc,h_1,p_1}\cdot
\big[\,\M_{0,(1^2),(2),(1^2)}^{V_1,V_0,V_2}(\F_1\cap U_0,2S;C_{pt^2},C_F,C_{pt^2})\,\big]
\cdot GT_{(1^2)}^{loc,h_2,p_2}
\end{align}
where the first equality follows from (\ref{pt-constraint}), Lemma~\ref{IL-scd}\,(c),
Remark~\ref{partition} and (\ref{connected-2}) and
the second equality follows from    Lemma~\ref{comp-ruled}.
Thus,
by Lemma~\ref{IL-scd}\,(c), (\ref{sum-scd-2}), (\ref{sum-scd-1}), Remark~\ref{lastremark} and Lemma~\ref{com:deg=1,2},
we have
\begin{align*}
GT_{(2)}^{loc,h,p}\,&=\,
(-1)^{p_1}2^{h_1}\,  GT_{(2)}^{loc,h_2,p_2} \, +\,
(-1)^{p_2}2^{h_2}\, GT_{(2)}^{loc,h_1,p_1}\\
&\ \ +\,
(-1)^p 2^{h-2}\, \big[\,\M_{0,(1^2),(2),(1^2)}^{V_1,V_0,V_2}(\F_1\cap U_0,2S;C_{pt^2},C_F,C_{pt^2})\,\big].
\end{align*}
When $(h_2,p_2)=(0,+)$, this equation shows
\begin{equation}\label{com-by-rel}
\big[\,\M_{0,(1^2),(2),(1^2)}^{V_1,V_0,V_2}(\F_1\cap U_0,2S;C_{pt^2},C_F,C_{pt^2})\,\big]
\,=\,-4\,GT_{(2)}^{loc,0,+}.
\end{equation}
This completes the proof.
\qed

\vskip 1cm
\setcounter{equation}{0}
\section{Spin Curve Degeneration and Sum Formula II}
\label{S:SC-Sum-2}
\bigskip

This section proves Theorem~B\,(b).
Let $h\geq 2$ or $(h,p)=(1,+)$ and
let ${\cal D}\to \Delta$ denote a family of curves over the unit disk $\Delta$ in $\cx$ whose fiber
over $\lam\ne 0$ is a smooth curve $D_\lam$ of genus $h$ and whose central fiber $D_0$ is a union of
two smooth components $\bar{D}$ and $E$ of genera $h-1$ and $0$, meeting at two points.
Fix a theta characteristic $\bar{N}$ on $\bar{D}$ with parity $p$. One can then find a line bundle ${\cal N}\to {\cal D}$
that restricts to a theta characteristic on $D_\lam$ with parity $p$,
to the theta characteristic $\bar{N}$ on $\bar{D}$  and to $\O(1)$ on $E$ (cf. pg 570 \cite{C}).
The projectivization $\P({\cal N}\oplus \O_{\cal D})$ gives a degeneration
$$
\P({\cal N}\oplus \O_{\cal D})\,\to\,{\cal D}\,\to\,\Delta
$$
such that (i) the general fiber $Z_\lam$ $(\lam\ne 0)$  is a ruled surface isomorphic to $\P_h=\P(N\oplus \O_D)$
where $(D,N)$ is a smooth spin curve of genus $h$ with parity $p$ and
(ii) the central fiber is the singular (ruled) surface
$$
\P_{h-1}\,\underset{V_1\sqcup\, V_2}{\textstyle \bigcup}\F_1\,\to\,D_0
$$
where $\P_{h-1}=\P(\bar{N}\oplus\O_{\bar{D}})$ and $V_1$ and $V_2$ are fibers over the nodes of $D_0$.
Note that the general fiber  $Z_\lam$ $(\lam\ne 0)$ is the symplectic fiber sum
\begin{equation*}\label{sym-sum-3}
\P_h\,=\,\P_{h-1}\,\underset{V_1\sqcup\, V_2}{\#}\,\F_1.
\end{equation*}

\bigskip

\non
{\bf Proof of Theorem~B\,(b)\,:}
The proof is also identical to those of Theorem~A and Theorem~B\,(a).
We only sketch the proof.
Fix a normal disk
$B\subset {\cal D}$ to $E$ at some point that is not a nodal point of $D_0$ and for each small $|\lam|$,
let $\tilde{V}_\lam$ be the fiber of $Z_\lam\to D_\lam$ over
the intersection point of $B$ and $D_\lam$. Choose small neighborhoods  $U$ and $W$ of the zero section
of $Z=\P({\cal N}\oplus \O_{\cal D})$ satisfying $\overline{W}\subset U$ and set
$$
W_\lam\,=\,W\,\cap\,Z_\lam\ \ \ \ \ \ \
\mbox{and}\ \ \ \ \ \ \
U_\lam\, =\,U\,\cap\, Z_\lam.
$$
The same arguments as in Lemma~\ref{IL-scd} and Lemma~\ref{key-scd} then show that
 the tautological section of $\pi^*{\cal N}$ over the total
space of $\pi:{\cal N}\to {\cal D}$ induces
an almost complex structure $J_{\sv}$ on $Z$ satisfying\,:
for $(J,\nu)$ sufficiently close to $(J_{\sv},0)$ and for small $|\lam|>0$
\begin{align}
&GT_{(2)}^{loc,h,p}\,=\,\big[\,\M_{\chi,(2)}^{\tilde{V}_\lam,*}(U_\lam,2S)\,\big]\ \ \ \ \ \
\mbox{and}\ \ \ \ \ \
GT_{s^1,s^2}^{loc,h-1,p}=\big[\,\M_{\chi_0,s^1,s^2}^{V_1,V_2,*}(\P_{h-1}\cap U_0)\,\big]
\label{IL-last}\\
&\ \ \ \ \ \ \ \ \ \ \ \ \ \ \ \ \ \ \ \ \ \ \ \ \ \ \ \ \
\M_{\chi,(2)}^{\tilde{V}_\lam,*}(U_\lam,2S) \setminus
\M_{\chi,(2)}^{\tilde{V}_\lam,*}(W_\lam,2S)=\emptyset^{\phantom{\displaystyle \sum}}
\label{key-3}
\end{align}
where $\chi=2-4h$ and $\chi_0=-4h+2\sum \ell(s^i)$.
For ordered sequences $s^i$ with $\deg(s^i)=2$ where $i=1,2$,
consider the evaluation map that records the intersection points with $V_1$ and $V_2$\,:
\begin{align*}
ev_{s^1,s^2, U_0}^* \,: \ &{\textstyle \bigcup}\
\M_{\chi_0,s^1,s^2}^{V_1,V_2,*}(\P_{h-1}\cap U_0,2S)\times
\M_{\chi_1,s^1,(2),s^2}^{V_1,V_0,V_2,*}(\F_1\cap U_0,2S)
\\
&\ \longrightarrow\
\Big(V_1^{\ell(s^1)}\times V_1^{\ell(s^1)}\Big )\times
\Big(V_2^{\ell(s^2)}\times V_2^{\ell(s^2)}\Big )^{\phantom{\displaystyle \cup}}
\end{align*}
where $V_0=\tilde{V}_0$, the union is over all
$ \chi_1+\chi_0 -2\ell(s^1)-2\ell(s^2)=2-4h.$
Let $\triangle_{s^i}$ be the diagonal of $V^{\ell(s^i)}\times V^{\ell(s^i)}$ where $i=1,2$.
Then we have
\begin{align}
&GT_{(2)}^{loc,h,p}\,=\,
[\M_{\chi,(2)}^{\tilde{V}_\lam,*}(U_\lam,2S)]\,=\,
\sum_{s^1,s^2}
\frac{|s^1||s^2|}{\ell(s^1)!\,\ell(s^2)!}\,
\big[\,\big(ev_{s^1,s^2,U_0}^*\big)^{-1}(\triangle_{s^1}\times\triangle_{s^2} )\,\big]
\notag \\
=&\
\sum_{m^1,m^2}\frac{|m^1||m^2|}{m^1!\,m^2!}\,\,
GT_{m^1,m^2}^{loc,h-1,p}\cdot
\big[\,\M_{0,m^1,(2),m^2}^{V_1,V_0,V_2}(\F_1\cap U_0,2S;C_{pt^{\ell(m^1)}},C_F,C_{pt^{\ell(m^2)}})\,\big]
\label{final}
\end{align}
where the first sum is over all ordered sequences $s^1$ and $s^2$
with $\deg(s^1)=\deg(s^2)=2$ and the second sum is over all partitions $m^1$ and $m^2$ of 2;
the first equality follows from (\ref{IL-last}), the second equality  from
(\ref{key-3}) and Theorem~10.1 of \cite{IP2} and the third equality
from (\ref{pt-constraint}), (\ref{IL-last}),
Remark~\ref{partition} and (\ref{connected-2}).
On the other hand, by Remark~\ref{lastremark} we have
$$
GT_{(1^2),(1^2)}^{loc,h-1,p}\,=\,2^2\,GT_{2}^{loc,h-1,p}\ \ \ \ \
\mbox{and}\ \ \ \ \
GT_{(1^2),(2)}^{loc,h-1,p}\,=\,2\,GT_{(2)}^{loc,h-1,p}.
$$
This together with (\ref{final}), Lemma~\ref{com:deg=1,2},
Lemma~\ref{comp-ruled} and (\ref{com-by-rel}) completes the proof. \qed

\vskip 1cm

\setcounter{equation}{0}
\section{Reduction to Genus Zero Spin Curve Invariants}
\label{S:MP}
\bigskip

As described in the Introduction,
Kiem and Li proved the Maulik-Pandharipande formulas (\ref{MPf}) by reducing
higher genus spin curve invariants to genus zero spin curve invariants.
The aim of this section is to show how their reduction follows from Theorem~A and Theorem~B.

\begin{prop}[\cite{KL1}]\label{P:Main}

\begin{itemize}
\item[]
\item[(a)]
${\displaystyle GT_1^{loc,h,p}\big(\,\prod_{i=1}^n\tau_{k_i}(F^*)\,\big)
\,=\, (-1)^p\,GT_1^{loc,0,+}\big(\,\prod_{i=1}^n\tau_{k_i}(F^*)\,\big)
}$
\item[(b)]
${\displaystyle GT_2^{loc,h,p}\big(\,\prod_{i=1}^n\tau_{k_i}(F^*)\,\big)
\,=\, (-1)^p \,2^h\,
GT_2^{loc,0,+}\big(\,\prod_{i=1}^n\tau_{k_i}(F^*)\,\big)
}$
\end{itemize}
\end{prop}

\pf
The sum formula (\ref{sum-bu}) for $d=1,2$ and Remark~\ref{lastremark} show that
\begin{align}
GT_1^{loc,h,p}\big(\,\prod_{i=1}^n\tau_{k_i}(F^*)\,\big)
&\,=\, GT_{1}^{loc,h,p}\cdot
GT_{(1),(1)}^{\F_0}\big(\,\prod_{i=1}^n\phi_i^{k_i}(F^*)\,\big)
\label{sum-1}  \\
GT_2^{loc,h,p}\big(\,\prod_{i=1}^n\tau_{k_i}(F^*)\,\big) &\,=\,
\frac{1}{2}\,GT_{2}^{loc,h,p}\cdot
GT_{(1,1),(1,1)}^{\F_0}\big(\,\prod_{i=1}^n\phi_i^{k_i}(F^*)\,\big)
\notag \\ &\   \,+\,
GT_{(2)}^{loc,h,p}\cdot
GT_{(2),(1,1)}^{\F_0}\big(\,\prod_{i=1}^n\phi_i^{k_i}(F^*)\,\big).
\label{sum-2}
\end{align}
Thus, Proposition~\ref{P:Main}\,(a) follows from Lemma~\ref{com:deg=1,2} and (\ref{sum-1}).
Similarly, by Lemma~\ref{com:deg=1,2} and (\ref{sum-2}), in order to prove Proposition~\ref{P:Main}\,(b),
we need to show
\begin{equation}\label{E:Main}
GT_{(2)}^{loc,h,p} \,= \,  (-1)^p \,2^h\, GT_{(2)}^{loc,0,+}.
\end{equation}
The sum formula (\ref{sum-deg-1}) for the case $(h_2,p_2)=(1,+)$ gives
\begin{equation}\label{sum-ind}
GT_{(2)}^{loc,h,p}\,=\,(-1)^p \,2^{h-1} GT_{(2)}^{loc,1,+}\,+\,2\, GT_{(2)}^{loc,h-1,p}\,-\,(-1)^p \,2^h GT_{(2)}^{loc,0,+}
\end{equation}
where $h\geq 2$. Applying the sum formula (\ref{sum-deg-1}) twice with $p_1=p_2=\pm 1$ and $h_1=h_2=1$ gives
$$
GT_{(2)}^{loc,1,+}\,=\,-GT_{(2)}^{loc,1,-}.
$$
This together with the sum formula (\ref{sum-deg-2}) for the case $(h,p)=(1,+)$ yields
\begin{equation}\label{IV-2}
GT_{(2)}^{loc,1,p}\,  =\, (-1)^p \,2\, GT_{(2)}^{loc,0,+}.
\end{equation}
Using induction on genus $h$ together with (\ref{sum-ind}) and (\ref{IV-2}) then shows (\ref{E:Main}). This
completes the proof. \qed

\medskip
\begin{rem} The proof of Proposition~\ref{P:Main}\,(b) by Kiem and Li (see Section 4 of \cite{KL1})
goes as follows\,: they first obtained a sum formula similar to  (\ref{sum-2}) using their sum formula and
then showed  (\ref{E:Main}) by calculating the local invariants
 $GT_2^{loc,h,p}(\tau(F^*))$ for all $h\geq 0$
using explicit algebro-geometric arguments.
\end{rem}

\vskip 1cm

\setcounter{equation}{0}
\section{Appendix}
\label{S:Appendix}
\bigskip

Let
$GW_{d,g}(\,\cdot\,)$ and $GW_{(1^d),(1^d),g}(\,\cdot\,)$ respectively denote
the absolute  GW invariants of $\P_0$ for the class $dS$
with genus $g$ and the relative GW invariants of $\P_0$ relative to
distinct fibers $V_1$ and $V_2$ of $\P_h$ with contact constraint $C_{[V_i]^d}$ with $V_i$
(we will omit the fibers $V_i$ and the contact constraints $C_{[V_i]^d}$ in notation).
Since local invariants of spin curve of genus $h=0$ are GW invariants of $\P_0$,
the lemma below shows the formula (\ref{E;step2R})  for the case when $h=0$ and $n\geq 3$.

\begin{lemma}\label{ap}
For $n\geq 3$, we have
\begin{equation}\label{E:ap}
GW_{d,g}\big(\,\prod_{i=1}^n \tau_{k_i}(F^*)\,\big)\,=\,
\frac{1}{(d!)^2}\,GW_{(1^d),(1^d),g}\big(\,\prod_{i=1}^n \phi_i^{k_i}(F^*)\,\big).
\end{equation}
\end{lemma}

\bigskip

The proof consists of two steps.

\bigskip
\non
{\bf Step 1 :}
We will relate the descendent classes for GW invariants of $\P_0$ to
the $\phi_i$ classes.
Following \cite{KM}, we set
$$
\tau_{s_i}\phi_i^{t_i}(F^*)\,=\,\psi_i^{s_i}st^*\phi_i^{t_i}\cup ev_i^*(F^*).
$$

\begin{lemma}\label{L:TRR}
Let $n\geq 3$. Then, for $s_j\geq 1$, we have
\begin{equation}\label{TRR}
\begin{array}{c}
{\displaystyle
GW_{d,g}\big(\,\prod_{i=1}^n \tau_{s_i}\phi_i^{t_i}(F^*)\,\big) \,=\,
GW_{d,g}\big(\,\prod_{i=1}^n \tau_{s_i-\delta_{ij}}\phi_i^{t_i+\delta_{ij}}(F^*)\,\big) } \\
{\displaystyle
-\,\sum_{0<k<d} \delta_{ks_j}\,(-1)^{k-1}\frac{1}{k!}\,
GW_{d-k,g}\big(\,\phi_j^{t_j}(F^*)\prod_{i\ne j} \tau_{s_i}\phi_i^{t_i}(F^*)\,\big).
}
\end{array}
\end{equation}
\end{lemma}

\pf
It follows from Theorem 1.1 of \cite{KM} that for $s_j\geq 1$
\begin{equation}\label{TRR1}
\begin{array}{c}
{\displaystyle
GW_{d,g}\big(\,\prod_{i=1}^n \tau_{s_i}\phi_i^{t_i}(F^*)\,\big) \,=\,
GW_{d,g}\big(\,\prod_{i=1}^n \tau_{s_i-\delta_{ij}}\phi_i^{t_i+\delta_{ij}}(F^*)\,\big) } \\
{\displaystyle
\,+\,\sum_{a,\,0<k\leq d}
GW_{k,0}\big(\,\tau_{s_j-1}(F^*)H^a\,\big)\,
GW_{d-k,g}\big(\,\phi_j^{t_j}(H_a)\prod_{i\ne j} \tau_{s_i}\phi_i^{t_i}(F^*)\,\big)
}
\end{array}
\end{equation}
where $\{H_a\}$ and $\{H^a\}$ are Poincar\'{e} dual basis of $H^*(\P_0)$.
Fix a basis  $\{1,S^*+F^*,F^*,\gamma^*\}$ and its dual basis
$\{\gamma^*,F^*,S^*,1\}$ of $H^*(\P_0)$ where  $\gamma^*$ be the Poincar\'{e} dual of the point class of $\P_0$.
Note that
all degree zero ($d=0)$ invariants in the righthand side of (\ref{TRR1}) vanish since $n\geq 3$;
no degree zero maps can pass through two distinct fibers.
Moreover, for any $d>0$ and $g$
$$
GW_{d,g}(\gamma^*\,\cdots\,)=0,\ \ \ \
GW_{d,g}(\,(S^*+F^*)\,\cdots\,)=0,\ \ \ \
GW_{d,g}(\,S^*\,\cdots\,)=-GW_{d,g}(\,F^*\,\cdots\,)
$$
where the first follows from $S^2=-1$, the second from $S(S+F)=0$ and the third from the second.
Consequently, (\ref{TRR1}) becomes
\begin{equation}\label{TRR2}
\begin{array}{c}
{\displaystyle
GW_{d,g}\big(\,\prod_{i=1}^n \tau_{s_i}\phi_i^{t_i}(F^*)\,\big) \,=\,
GW_{d,g}\big(\,\prod_{i=1}^n \tau_{s_i-\delta_{ij}}\phi_i^{t_i+\delta_{ij}}(F^*)\,\big) } \\
{\displaystyle
-\,\sum_{0<k<d} GW_{k,0}\big(\,\tau_{s_j-1}(F^*)F^*\,\big)\,
GW_{d-k,g}\big(\,\phi_j^{t_j}(F^*)\prod_{i\ne j} \tau_{s_i}\phi_i^{t_i}(F^*)\,\big).
}
\end{array}
\end{equation}
If $k\ne s_j$ then $GW_{k,0}(\tau_{s_j-1}(F^*)F^*)=0$ by dimension count.
So, it remains to show
\begin{equation}\label{lasteq}
GW_{k,0}\big(\tau_{k-1}(F^*)F^*)\,=\,(-1)^{k-1}/k!
\end{equation}
When $k=s_j$, the generalized Divisor Axiom (cf. Lemma 1.4 of \cite{KM}) and
(\ref{TRR2}) (applied to $GW_{k,0}(\tau_{k-1}(F^*)F^*F^*)$)
together with the facts $\phi_1=0$ on $\CM_{0,3}$ and $GW_{1,0}(F^*F^*)=1$ give
$$
GW_{k,0}\big(\,\tau_{k-1}(F^*)F^*\,\big)
\,=\,\frac1k\,GW_{k,0}\big(\,\tau_{k-1}(F^*)F^*F^*\,\big)
\,=\,-\frac1k\,GW_{k-1,0}\big(\,\tau_{k-2}(F^*)F^*\,\big)
$$
By induction, this shows (\ref{lasteq}) that completes the proof. \qed

\bigskip
\non
{\bf Step 2 :}
We first show a formula for relative invariants that is analogous to (\ref{TRR}) and then give a proof of
Lemma~\ref{ap}.  Recall that  for the
forgetful map $\pi_{\ell}:\CM_{g,n+\ell}\to\CM_{g,n}$ that forgets the last $\ell$ marked points and for
$1\leq i\leq n$ we have
\begin{equation}\label{pullbackclassbyxi-a}
\pi_{\ell}^*\,\phi_i\,=\,\phi_i-\sum \delta_{\{i\}\cup I}
\end{equation}
where the sum is over all $I\subset \{n+1,\ldots,n+\ell\}$ with $I\ne \emptyset$.
For simplicity, we will write $\pi_{\ell}$ simply as $\pi$ when $\ell$ is even.

\begin{lemma}\label{L:TRR-rel}
Let $n\geq 3$. Then, for $s_j\geq 1$ we have
\begin{equation*}
\begin{array}{c}
{\displaystyle
GW_{(1^d),(1^d),g}\big(\,\prod_{i=1}^n \phi_i^{s_i}\pi^*\phi_i^{t_i}(F^*)\,\big) \,=\,
GW_{(1^d),(1^d),g}\big(\,\prod_{i=1}^n \phi_i^{s_i-\delta_{ij}}\pi^*\phi_i^{t_i+\delta_{ij}}(F^*)\,\big) } \\
{\displaystyle
\,-\,\sum_{0<k< d}\, \delta_{ks_j}\,(-1)^{k-1}\,k!\,
\left(\begin{array}{c} d\\ k \end{array} \right)^2\,
GW_{(1^{d-k}),(1^{d-k}),g}\big(\,\pi^*\phi_j^{t_j}(F^*)\prod_{i\ne j} \phi_i^{s_i}\pi^*\phi_i^{t_i}(F^*)\,\big)
}
\end{array}
\end{equation*}
\end{lemma}

\pf
Without loss of generality, we may assume $j=1$.
For the forgetful map $\pi=\pi_{2d}$, let $\delta_{\{1\}\cup I}$ be a class as in (\ref{pullbackclassbyxi-a}) and
denote by  $\CM(\delta_{\{1\}\cup I})$ the boundary stratum of $\CM_{g,n+2d}$
whose fundamental class is Poincar\'{e} dual to  $\delta_{\{1\}\cup I}$.
Then  for $m=|I|$ there is a gluing map
$$
\eta_{\si}\,:\, \CM_{0,m+2}\times \CM_{g,n+2d-m} \,\to\, \CM_{g,n+2d}
$$
whose image is
$\CM(\delta_{\{1\}\cup I})$. This gluing map is obtained by
identifying the second marked point of the first component  with the first marked point of the second component. We have
\begin{equation}\label{s-phi-eta}
\eta_{\si}^*(\phi_1)\,=\,\phi_1\otimes 1
\ \ \ \ \ \mbox{and}\ \ \ \ \
\eta_{\si}^*\circ \pi^*(\phi_1)=1\otimes \pi^*_{2d-m}(\phi_1)
\end{equation}
where the first equality follows from Lemma~\ref{L:AC}\,(b) and the second from the fact
that  under the composition map $\pi\circ \eta_{\si}$
the first component collapses to a point.

Choose two distinct fibers $V_1$ and $V_2$ of $\P_0$ and, for simplicity, we set
$$
\M^V\,=\,\M^{V_1,V_2}_{g,n,(1^d),(1^d)}(\P_0,dS)\ \ \ \ \
\mbox{and}\ \ \ \ \
\Phi\,=\,\phi_1^{s_1-1}\pi^*\phi_1^{t_1}\prod_{i>1} \phi_i^{s_i}\pi^*\phi_i^{t_i}
$$
where $V=V_1\sqcup V_2$.
Let $G$ be a geometric representative of the Poincar\'{e} dual of the pull-back class $\eta_{\si}^*\Phi$.
One can then choose a (smooth) family of geometric representatives $G_t$ of the Poincar\'{e} dual of the class
$\delta_{\{1\}\cup I}\cup \Phi$ with $G_0=\eta_{\si}(G)$.
Let $B$ be a product of $n$ distinct generic fibers $B_i$ of $\P_0$ each of which is  disjoint with $V$.

Suppose $\M^V\cap B\cap G_t\ne \emptyset$ for all small $t$.
Then, by the Gromov Compactness Theorem, after passing to subsequences,
as $t\to 0$ every sequence $f_t\in \M^V\cap B\cap G_t$ converges to
$$
(f,C)\ \in \ C\M^V\,\cap\,B\,\cap\,\eta_{\si}(G)
$$
where   $C\M^V$ is the closure of $\M^V$ in $\CM_{g,n+2d}(\P_0,dS)$.
The closure $C\M^V$ has a stratification in which each stratum consisting of maps with domains with more than one node
has dimension at least $4$ less than $2\deg(\Phi)+2+2n$ (cf. Lemma 7.6 of \cite{IP2}).
Thus the domain $C \in \CM(\delta_I)$ of $f$ has one node by dimension count.
The limit map $f$ splits as $f=(f_1,f_2)$ such that each $f_i$ is $V$-regular unless it represents the trivial homology class.
In our case, both $f_1$ and $f_2$ are $V$-regular maps since
the image of $f_1$ passes through $V$ and $B_1$, and the image of $f_2$ passes through $(n-1)>2$ distinct fibers
$B_i$ where $2\leq i\leq n-1$. For some $0<k<d$, we have
\begin{itemize}
\item
$f_1$ (resp. $f_2$) has contact vector $(1^k)$ (resp. $(1^{d-k})$) with both $V_1$ and $V_2$, and hence
\item
$f\in ev_{\si}^{-1}(\triangle)\cap B \cap \eta(G)$
(under the natural inclusion $ev_{\si}^{-1}(\triangle) \hookrightarrow \M^V$)
\end{itemize}
where $\triangle$ is the diagonal of $\P_0\times \P_0$ and $ev_{\si}$ is the evaluation map
$$
ev_{\si}=ev_{2}\times ev_{1}\,:\,
\M^{V_1,V_2}_{0,2,(1^k),(1^k)}(\P_0,kS)\,\times
\M^{V_1,V_2}_{g,n,(1^{d-k}),(1^{d-k})}(\P_0,(d-k)S)\ \mapsto\ \P_0\times\P_0.
$$
On the other hand, the condition of contact order $(1^d)$ with $V$ is an open condition.
The  Gluing Theorem of \cite{RT2} thus implies that
for each small $t$ one can uniquely smooth $f$ (at a node) to produce a $V$-regular map in $\M^V\cap G_t$.
Consequently, we have
\begin{equation}\label{key-composition}
\big[\,ev_{\si}^{-1}(\triangle)\,]\,\cap\,\eta_{\si}^*(\Phi) \otimes (F^*)^n \,=\,
\big[\,\M^V\,]\,\cap\ (\delta_{\{1\}\cup I} \cup\Phi)\otimes (F^*)^n.
\end{equation}

It follows from (\ref{pullbackclassbyxi-1}), (\ref{s-phi-eta}), (\ref{key-composition}),
the splitting of the diagonal $\triangle$ as in the proof of Lemma~\ref{L:TRR} (see paragraph above (\ref{TRR2}))
that
\begin{equation}\label{TRR-rel}
\begin{array}{c}
{\displaystyle
GW_{(1^d),(1^d),g}\big(\,\prod_{i=1}^n \phi_i^{s_i}\pi^*\phi_i^{t_i}(F^*)\,\big) \,=\,
GW_{(1^d),(1^d),g}\big(\,\prod_{i=1}^n \phi_i^{s_i-\delta_{i1}}\pi^*\phi_i^{t_i+\delta_{i1}}(F^*)\,\big) } \\
{\displaystyle
\,-\,\sum_{0<k< d}\,
\Big(\begin{array}{c} d \\ k \end{array} \Big)^2\,
GW_{(1^k),(1^k), 0}\big(\,\phi_1^{s_1-1}(F^*)F^*\,)\cdot
GW_{(1^{d-k}),(1^{d-k}),g}\big(\,\pi^*\phi_1^{t_1}(F^*)\prod_{i>1} \phi_i^{s_i}\pi^*\phi_i^{t_i}(F^*)\,\big)
}
\end{array}
\end{equation}
where the factor $\Big(\!\begin{array}{c} d \\ k \end{array} \!\Big)^2$ reflects the fact that  the classes
$\delta_{\{1\}\cup I}$ in (\ref{key-composition}) are obtained by choosing $k$ contact points  with $V_1$ and
$k$ contact points with $V_2$.
Observe that by dimension if $k\ne s_1$ then $GW_{(1^k),(1^k),0}\big(\phi_1^{s_1-1}(F^*)F^*)=0.$
Thus, it remains to show\,:
\begin{equation}\label{interesting}
GW_{(1^k),(1^k),0}\big(\,\phi^{k-1}(F^*)F^*\,)\,=\,(-1)^{k-1}k!
\end{equation}
The sum formula (\ref{oursum}) for $h=0$ gives
\begin{align}\label{last-com1}
GW_{(1^k),(1^k),0}\big(\,\phi^{k-1}(F^*)F^*F^*\,\big)\,&=\,
\frac{1}{k!}\,GW_{(1^k),(1^k),0}\big(\,\phi^{k-1}(F^*)F^*\,\big)\cdot
GT_{(1^k),(1^k)}^{\F_0}\big(\,F^*\,\big)
\notag \\
&=\,
k\,GW_{(1^k),(1^k),0}\big(\,\phi^{k-1}(F^*)F^*\,\big)
\end{align}
where  the second equality follows from (\ref{com-F0}) and the Divisor Axiom.
On the other hand, together with the facts $\phi=0$ on $\CM_{0,3}$ and $GW_{(1),(1),0}(F^*F^*)=1$,
the formula (\ref{TRR-rel}) shows
\begin{equation}\label{last-com2}
GW_{(1^k),(1^k),0}\big(\,\phi^{k-1}(F^*)F^*F^*\,\big)\,=\,
-k^2\,GW_{(1^{k-1}),(1^{k-1}),0}\big(\,\phi^{k-2}(F^*)F^*\,\big).
\end{equation}
By induction, (\ref{last-com1}) and (\ref{last-com2}) thus imply (\ref{interesting}). This completes the proof. \qed

\bigskip
\non
{\bf Proof of Lemma~\ref{ap} :} It suffice to show that for $n\geq 3$
\begin{equation}\label{induction}
GW_{d,g}\big(\,\prod_{i=1}^n \tau_{s_i}\phi_i^{t_i}(F^*)\,\big) \,=\,
\frac{1}{(d!)^2}\,GW_{(1^d),(1^d),g}\big(\,\prod_{i=1}^n \phi_i^{s_i}\pi^*\phi_i^{t_i}(F^*)\,\big).
\end{equation}
When $\sum s_i=0$, (\ref{induction})  follows from Lemma~\ref{step1h=0}.
Suppose that (\ref{induction}) holds  for any $d$, $g$ and $n\geq 3$ whenever $\sum s_i<\ell$.
Then, Lemma~\ref{L:TRR} and Lemma~\ref{L:TRR-rel} show that
(\ref{induction}) also holds when $\sum s_i=\ell$.
Therefore, (\ref{induction}) follows from induction on the sum $\sum s_i$.
\qed

\end{document}